\documentclass[11pt]{article}
\usepackage{amssymb, amsthm, amsmath, amscd}
\newtheorem{thm}{Theorem}[section]
\newtheorem{lem}[thm]{Lemma}
\newtheorem{prop}[thm]{Proposition}
\newtheorem{cor}[thm]{Corollary}
\newtheorem{NN}[thm]{}
\theoremstyle{definition}\newtheorem{df}[thm]{Definition}
\theoremstyle{definition}\newtheorem{rem}[thm]{Remark}
\theoremstyle{definition}

\renewcommand{\phi}{\varphi}

\newcommand{\N}{\mathbb{N}}
\newcommand{\Z}{\mathbb{Z}}
\newcommand{\Q}{\mathbb{Q}}
\newcommand{\R}{\mathbb{R}}

\newcommand{\T}{\mathbb{T}}

\newcommand{\id}{\operatorname{id}}

\newcommand{\morp}{contractive completely positive linear map}

\newcommand{\hm}{homomorphism}
\newcommand{\dt}{\delta}
\newcommand{\ep}{\epsilon}
\newcommand{\andeqn}{\,\,\,{\rm and}\,\,\,}
\newcommand{\rforal}{\,\,\,{\rm for\,\,\,all}\,\,\,}
\newcommand{\CA}{$C^*$-algebra}
\newcommand{\SCA}{$C^*$-subalgebra}

\newcommand{\af}{{\alpha}}
\newcommand{\bt}{{\beta}}

\newcommand{\D}{\mathbb D}
\newcommand{\beq}{\begin{eqnarray}}
\newcommand{\eneq}{\end{eqnarray}}
\newcommand{\tforal}{\,\,\,\text{for\,\,\,all}\,\,\,}

\title{AF-embedding  of crossed products of AH-algebras by $\Z$ and
asymptotic AF-embedding}
\author{Huaxin Lin
 }
\date{}

\begin{document}

\maketitle

\begin{abstract}
Let $A$ be a unital AH-algebra and let $\af\in Aut(A)$ be an
automorphism.  A necessary condition for $A\rtimes_{\af}\Z$ being
embedded into a unital simple AF-algebra is the existence of a
faithful tracial state. If in addition, there is an automorphism
$\kappa$ with $\kappa_{*1}=-{\rm id}_{K_1(A)}$ such that $\af\circ
\kappa$ and $\kappa\circ \af$ are asymptotically unitarily
equivalent, then $A\rtimes_{\af}\Z$ can be embedded into a unital
simple AF-algebra.  Consequently,
in the case that  $A$ is a unital AH-algebra (not necessarily
simple) with torsion $K_1(A),$ $A\rtimes_{\af}\Z$ can be embedded
into a unital simple AF-algebra if and only if $A$ admits a
faithful $\af$-invariant tracial state. We also show that if $A$
is a unital A$\T$-algebra then $A\rtimes_{\af}\Z$ can be embedded
into a unital simple AF-algebra if and only if $A$ admits a
faithful $\af$-invariant tracial state. Consequently, for any unital simple A$\T$-algebra
$A,$ $A\rtimes_{\af}\Z$ can always be embedded into a unital simple AF-algebra.

If $X$ is a compact metric space and $\Lambda: \Z^2\to Aut(C(X))$
is a \hm\, then $C(X)\rtimes_{\Lambda}\Z^2$ can be asymptotically
embedded into a unital simple AF-algebra provided that $X$ admits
a strictly positive $\Lambda$-invariant probability measure.
Consequently $C(X)\rtimes_{\Lambda}\Z^2$ is quasidiagonal if $X$
admits a strictly positive $\Lambda$-invariant Borel probability
measure.

\end{abstract}

\section{Introduction}
Quasidiagonality for crossed product \CA s were studied by Dan Voiculescu (\cite{V2}), \cite{V3} and \cite{V4}). Quasidiagonality in \CA s has been studied for a long time in many different point of views
(see \cite{Bn2} for more information, also \cite{DHS}, \cite{Th}, \cite{Hd}, \cite{Ros}, \cite{Sa}, \cite{Z}, \cite{BK1}, \cite{BK2}, \cite{BlD}, \cite{ELP}, \cite{D1} and \cite{D2}, to name a few).
Let $X$ be a
compact metric space and $\af$ be a homeomorphism on $X.$  It was
proved by Pimsner (\cite{Pi}) that $C(X)\rtimes_{\af}\Z$ is
quasidiagonal  if and only if $\af$ is pseudo-non-wondering, and
if and only if $C(X)\rtimes_{\af}\Z$ can be embedded into an
AF-algebra.
Recent applications of \CA\, theory to dynamical
systems further suggest that quasidiagonality and
$AF$-embeddability  of crossed products $A\rtimes_{\af}\Z$ become
increasingly important.

 Let $A$ be an AF-algebra and $\af\in Aut(A)$ be an automorphism.
Nate Brown (\cite{Bn1}) proved that $A\rtimes_{\af}\Z$ can be
embedded into an AF-algebra if and only if $A\rtimes_{\af}\Z$ is
quasidiagonal. More importantly, a $K$-theoretical necessary and sufficient
condition for $A\rtimes_{\af}\Z$ being quasidiagonal is also
given there. When $A$ is a unital simple $A\T$-algebra of real rank zero, Matui
(\cite{M1}) showed that $A\rtimes_{\af}\Z$ can always be embedded
into an AF-algebra. On the other hand, Voiculescu's question when
$C(X)\rtimes_{\Lambda}\Z^2$ can be embedded into an AF-algebra
 remains unanswered.  It was not even known when
 $C(X)\rtimes_{\Lambda}\Z^2$ is quasidiagonal in general. It was pointed out by Voiculescu ( \cite{V4}) that
 a unital quasidiagonal \CA\, always has a tracial state.
 One consequence of the results here is that if $C(X)\rtimes_{\af}\Z^2$ admits a faithful tracial state
 then it is quasidiagonal (\ref{Cqd}).

We are interested in the problem when \CA s can be embedded into
unital simple AF-algebras. Note that not every AF-algebra can be
embedded into a unital simple AF-algebra. In fact, it is easy to show
that an AF-algebra can be embedded into a unital simple AF-algebra
if and only if it admits a faithful tracial state. We prove  a
more general result. Recall that a \CA\, $A$ is AH,
if $A$ is an inductive limit of finite direct sums of \CA s $A_n$
with the form $P_{(n,i)} M_{k(n,i))}(C(X_{n,i}))P_{(n,i)},$ where
$X_{n,i}$ is a finite CW complex and $P_{(n,i)}\in
M_{k(n,i)}(C(X_{n,i}))$ is a projection. We show that a unital
AH-algebra can be embedded into a unital simple AF-algebra if and
only if it admits a faithful tracial state (\ref{embAH}).

Let $A$ be a unital AH-algebra and let $\af\in Aut(A)$ be an
automorphism. We study the problem when $A\rtimes_{\af}\Z$ can be
embedded into a unital simple AF-algebra. A necessary condition is
that $A$ admits a faithful $\af$-invariant tracial state. We
provide a sufficient condition for $A\rtimes_{\af}\Z$ being
embedded into a unital simple AF-algebra. We prove this by
applying our recent result of so-called Basic Homotopy Lemma in
\cite{Lnhomp}. We show that, when there exists an automorphism
$\kappa$ with $\kappa_{*1}=-{\rm id}_{K_1(A)}$ such that
$\kappa\circ \af$ and $\af\circ \kappa$ are asymptotically
unitarily equivalent, $A\rtimes_{\af}\Z$ can be embedded into a
unital simple AF-algebra if and only if $A$ admits a faithful
$\af$-invariant tracial state (\ref{AHembT1}). It should be noted
that such $\kappa$ do exist in many known cases.
It also shows that, in the case that  $A$ is a unital AH-algebra
(not necessarily simple) with torsion $K_1(A),$ $A\rtimes_{\af}\Z$
can be embedded into a unital simple AF-algebra if and only if $A$
admits a faithful $\af$-invariant tracial state. It should be
noted that a unital AH-algebra $A$ with zero $K_1(A)$ may not be
AF-algebra. Even in the simple case, a unital simple AH-algebra
with zero $K_1(A)$ and with torsion free $K_0(A)$ may not be
AF-algebra. There are simple AH-algebras with stable rank higher
than one (\cite{Vj}). On the other hand, we also show that when
$A$ is a unital A$\T$-algebra $A\rtimes_{\af}\Z$ can be embedded
into a unital simple AF-algebra if and only if $A$ admits a
faithful $\af$-invariant tracial state. In particular, when $A$ is a unital simple A$\T$-algebra,
$A\rtimes_{\af}\Z$ can always be embedded into a unital simple AF-algebra, which generalizes a result
of Matui (\cite{M1}).

 A \CA\,
$A$ is said to be asymptotically embedded into a \CA\, $B,$ if
there exists a sequence of \morp s $\{\phi_n\}: A\to B$ such that
$\lim_{n\to\infty}\|\phi_n(ab)-\phi_n(a)\phi_n(b)\|=0$ for all $a,
b\in A$ and $\lim_{n\to\infty}\|\phi_n(a)\|=\|a\|.$ Among other
things, we show that $C(X)\rtimes_{\Lambda}\Z^2$ can be
asymptotically embedded (in a stronger sense) into a unital simple
AF-algebra if (and only if) $X$ admits a strictly positive
$\Lambda$-invariant Borel probability measure (\ref{TembX}).
Consequently, in the case that $X$ admits a strictly positive
$\Lambda$-invariant Borel probability measure,
$C(X)\rtimes_{\Lambda}\Z^2$ is quasidiagonal.

\vspace{0.2in}

{\bf Acknowledgments} The preliminary work of this paper was done
when the author was in East China Normal University in the summer
2006 which was partially supported by Shanghai Academic Priority Disciplines.
This work is also partially supported by a grant from
National Science Foundation.

\section{Preliminaries}

\begin{NN}
{\rm

 Let $A$ be a unital \CA. Denote by $U(A)$ the unitary group of $A$ and denote by $U_0(A)$ the path connected
 component of $U(A)$ containing the identity. Denote by $Aut(A)$ the
 group of automorphisms on $A.$ If $u\in U(A),$ denote by ${\rm
 ad}\, u$ the inner automorphism defined by ${\rm ad}\,
 u(a)=u^*au$ for all $a\in A.$
 \vspace{0.1in}

Let $A$ be a unital \CA\, and let $\af\in Aut(A).$ Denote by
$A\rtimes_\af\Z$ the crossed product. In this paper, we will fix a
unitary and denote it by $u_\af\in A\rtimes_\af\Z$ for which ${\rm
ad}\, u_\af(a)=\af(a)$ for all $a\in A.$
 \vspace{0.1in}

 Let $A$ be a stably finite \CA. Denote by $T(A)$ the tracial state space of $A$
 and by $Aff(T(A))$ the normed space of all real affine continuous functions on $T(A).$
 Denote by $\rho_A: K_0(A)\to Aff(T(A))$ the positive \hm\, induced by
 $\rho_A([e])=(\tau\otimes T_k)(e),$ where $e$ is a projection in $A\otimes M_k$ and $T_k$
 is the standard trace on $M_k,$ $k=1,2,....$
 \vspace{0.1in}

Let $A$ and $B$ be two \CA s and let $h: A\to B$ be a \hm. We use
$h_{*i}$ for the induced \hm\, from $K_i(A)$ into $K_i(B),$
$i=0,1.$
 \vspace{0.1in}

 A \hm\, $h$ is said to be {\it full}, if $h(a)$ does not belong to any non-trivial (two-sided closed) ideal
 of $B$ for any $a\not=0.$
 \vspace{0.1in}

 Let $\{B_n\}$ be a sequence of \CA s. Denote by $l^{\infty}(\{B_n\})$ the $C^*$-algebraic
product of $\{B_n\}$ and $c_0(\{B_n\})$ the $C^*$-algebraic direct
sum. We also use $q_{\infty}(\{B_n\})$ for the quotient
$l^{\infty}(\{B_n\})/c_0(\{B_n\}).$
 \vspace{0.1in}

 A measure $\mu$ on $X$ is said to be {\it strictly positive} measure,
if $\mu(O)>0$ for any non-empty open subset $O\subset X.$

}
\end{NN}


\begin{df}\label{dsinj}
{\rm Let $A$ and $B$ be two \CA s. Suppose that $L_1, L_2: A\to B$
are  two maps, $\ep>0$ and ${\cal F}\subset A,$ we write
$$
L_1\approx_{\ep} L_2 \,\,\,\text{on}\,\,\, {\cal F},
$$
if
$$
\|L_1(a)-L_2(a)\|<\ep\rforal a\in {\cal F}.
$$

 Suppose that $L: A\to B$ is a
\morp. Let $\dt>0$ and ${\cal G}\subset A$ be a (finite) subset.
We say that $L$ is $\dt$-${\cal G}$-multiplicative if
$$
\|L(ab)-L(a)L(b)\|<\dt\tforal a, b\in {\cal G}.
$$
Let $\sigma>0.$ We say that $L$ is  $\sigma$-${\cal G}$-injective
if
$$
\|L(a)\|\ge \sigma\|a\|\tforal a\in {\cal G}.
$$

Let ${\cal P}\subset \underline{K}(A)$ be a finite subset. A
\morp\, $L: A\to B$ is said to be $({\cal G}, {\cal P}, \dt)$ map
if $L$ is $\dt$-${\cal G}$-multiplicative, $[L]_{\cal P}$ is well
defined (see 2.2 of \cite{Lnhomp}, 6.1.1 of \cite{Lnbk} and see
also 2.4 of \cite{D1}).

}
\end{df}

\begin{df}\label{dasmeb}
{\rm Let $A$ and $B$ be two \CA s. Suppose that $\{\phi_n\}$ is a
sequence of \morp s from $A$ to $B.$ We say $\{\phi_n\}$ is an
{\it asymptotic embedding}, if
$$
\lim_{n\to\infty}\|\phi_n(a)\phi_n(b)-\phi_n(ab)\|=0 \tforal a,
b\in A\andeqn
$$
$$
\lim_{n\to\infty}\|\phi_n(a)\|=\|a\|\tforal a\in A.
$$
}

\end{df}

\begin{df}\label{Tfull}

{\rm Let $A$ and $B$ be two unital \CA s. Let $T: A_+\setminus
\{0\}\to \N\times R_+$ be a map and ${\cal G}\subset A_+$ be a
subset.  Suppose that $\phi: A\to B$ is a \morp. We say $\phi$ is
$T$-${\cal G}$-full, if there are $x_1,x_2,...,x_{T(a)_{\N}}\in A$
such that
$$
\max\{\|x_i\|: 1\le i\le T(a)_{\N}\}\le T(a)_{\R}\andeqn
\sum_{i=1}^{T(a)_\N}x_i^*\phi(a)x_i=1_B
$$
for all $a\in {\cal G}.$

 Suppose that
$\{\phi_n\}$ is a sequence of \morp s from $A$ to $B$ such that
$$
\lim_{n\to\infty}\|\phi_n(ab)-\phi_n(a)\phi_n(a)\|=0 \tforal a,
b\in A.
$$
We say $\{\phi_n\}$ is a $T$-full asymptotic embedding if for any
finite subset ${\cal G}\subset A_+,$ there exists an integer $N>0$
such that $\phi_n$ is $T$-${\cal G}$-full for all $n\ge N.$  From \ref{fullP1},
it is clear that such $\{\phi_n\}$ is indeed an asymptotic embedding.

}
\end{df}

\begin{df}
{\rm  Denote by ${\cal N}$ the class of separable amenable \CA s which satisfy the Universal Coefficient Theorem.

}
\end{df}

\begin{df}\label{DU}
{\rm Denote by ${\cal U}$ throughout this paper the universal
UHF-algebra ${\cal U}=\otimes_{n\ge1}M_n.$

Let $\{e_{i,j}^{(n)}\}$ be the canonical matrix units for $M_n.$
Let $u_n\in M_n$ be the unitary matrix such that ${\rm ad}\,
u_n(e_{i,i}^{(n)})=e_{i+1,i+1}^{(n)}$ (modulo $n$). Let
$\sigma=\otimes_{n\ge 1} {\rm ad}\, u_n \in Aut({\cal U})$ be the
shift (see for example Example 2.2 of \cite{Bn1}). A fact that we
will use in this paper is the following  property that $\sigma$
has: For any integer $k>0,$ any $\ep>0$ and any finite subset
${\cal F}\subset {\cal U},$  there exist mutually orthogonal
projections $e_1,e_2,...,e_k\in {\cal U}$ such that $
\sum_{i=1}^ke_i=1_{\cal U},$ $ \|xe_i-e_ix\|<\ep$ for all $ x\in
{\cal F}$ and $\sigma(e_i)=e_{i+1}, i=1,2,..., k\,\,\, {\rm
(}e_{k+1}=e_1{\rm )}.$

Let $A$ be a unital \CA\, and let $\af\in Aut(A).$ In this paper, we say that
$A$ has the {\it cyclic Rokhlin property}, if for any integer $k>0,$ any $\ep>0$ and any finite
subset ${\cal F}\subset A,$  there exist mutually
orthogonal projections $e_1,e_2,...,e_k\in A$ such that
\beq
&&\sum_{i=1}^ke_i=1_{A},\\
&& \|xe_i-e_ix\|<\ep\,\rforal x\in {\cal F}\andeqn\\
&&\|\af(e_i)-e_{i+1}\|<\ep,\,\,\, i=1,2,..., k\,\,\, {\rm (}e_{k+1}=e_1{\rm )}.
\eneq

Suppose that $\af$ has the above cyclic Rokhlin property and $w\in U(A). $  Then
${\rm ad}\, w\circ \af$ also has the cyclic Rokhlin property (a fact communicated to us by
N. C. Phillips). Given $\ep>0$ and finite subset ${\cal F}\subset A$ and an integer $k>0,$  defined
${\cal F}_1=\{w\}\cup{\cal F}.$
Since $\af$ has the cyclic Rokhlin property, there are mutually orthogonal projections
$\{e_1,e_2,...,e_k\}$ such that
$\sum_{i=1}^k e_i=1_A,$ $\|\af(e_i)-e_{i+1}\|<\ep/2$ and $\|xe_i-e_ix\|<\ep/2$ for all $x\in {\cal F}_1.$
Then
\beq\label{DU-2}
\|{\rm ad}\, w\circ \af(e_i)-e_{i+1}\| &\le &\|{\rm ad}\, w\circ \af(e_i)-{\rm ad}\, w(e_{i+1})\|
+\|{\rm ad}\, w(e_{i+1})-e_{i+1}\|\\
&< &\ep/2+\ep/2=\ep.
\eneq
It follows that ${\rm ad}\, w\circ \af$ has the cyclic Rokhlin property.

}
\end{df}

\section{ Two Uniqueness Theorems}

The following is known.

\begin{lem}\label{inj1}
Let $A$ and $B$ be two \CA s. Suppose that $\phi_n: A\to B$ is a
sequence of \morp\, such that
$$
\lim_{n\to\infty}\|\phi_n(a)\phi_n(b)-\phi_n(ab)\|=0\tforal a,
b\in A.
$$
Then the following are equivalent:

(i) The sequence $\{\phi_n\}$ is an asymptotic embedding;

(ii) There exists $\sigma>0$ and a dense and  increasing sequence
of finite subset ${\cal G}_n\subset A$ such that each $\phi_n$ is
$\sigma$-${\cal G}_n$-injective;

(iii) The induced the map: $\Phi: A\to q_{\infty}(\{B\})$ defined
by $\pi\circ \{\phi_n\},$ where $\pi: l^{\infty}(\{B\})\to
l^{\infty}(\{B\})\to q_{\infty}(\{B\})$ is the quotient map, is a
monomorphism.

\end{lem}

\begin{proof}
It is clear that (i) implies (ii).

If (ii) holds, then $\pi\circ (\{\phi_n\}): A\to
q_{\infty}(\{B_n\})$ is a \hm. Moreover,
$$
\|\pi\circ (\{\phi_n\})(a)\|\ge \sigma\|a\|\tforal a\in A.
$$
It follows that $\pi\circ (\{\phi_n\})$ is a monomorphism.

If (iii) holds, then
$$
\|\pi\circ (\{\phi_n\})(a)\|=\|a\|\tforal a\in A.
$$
It follows that
$$
\lim_{n\to\infty}\|\phi_n(a)\|=\|a\|\tforal a\in A.
$$
\end{proof}

\begin{lem}\label{Ldense}
Let $A$ be a unital separable amenable  \CA\, and let $T:
A_+\setminus \{0\}\to (\N, \R_+)$ be a map. Then there exists an
increasing sequence of finite subsets $\{{\cal G}_n\}\subset A_+$
with its union being dense in $A_+$ satisfying the following: If
$B_n$ is a sequence of  unital \CA\, and if $\phi_n: A \to B_n$ is
a \morp\, such that
$$
\lim_{n\to\infty}\|\phi_n(ab)-\phi_n(a)\phi_n(b)\|=0
$$
for all $a\in A$ and each $\phi_k$ is $T$-${\cal G}_n$ full, then
the induced map $\Phi=\pi\circ \{\phi_n\}$ from $A$ to
$q_{\infty}(\{B_n\})$ is a full map, where $\pi:
\prod_{n=1}^{\infty}(\{B_n\}\to q_{\infty}(\{B_n\})$ is the
quotient map.
\end{lem}

\begin{proof}
First we note that $\Phi: A\to q_{\infty}(\{B_n\})$ is a \hm.

Let $\{{\cal G}_n'\}$ be an increasing sequence of finite subsets
of $A_+$ such that $\cup_{n=1}^{\infty} {\cal G}_n'$ is dense in
$A_+.$  For each integer $n>0,$ define
$f_{1/n}^{1-1/n}$ to be a positive continuous function on $[0,
\infty)$ such that $0\le f_{1/n}^{1-1/n}(t)\le 1,$
$f_{1/n}^{1-1/n}(t)=1$ if $t\ge 1-1/n,$ $f_{1/n}^{1-1/n}(t)=0$ if
$t\in [0, 1/n]$ and $f_{1/n}^{1-1/n}(t)={t-1/n\over{1-2/n}}$ for
$t\in [1/n, 1-1/n].$

Let ${\cal G}_n''=\{f_{1/n}^{1-1/n}(g): g\in A_+\cap {\cal
G}_n'\}.$ Define ${\cal G}_n={\cal G}_n'\cup {\cal G}_n''.$ Now
suppose that $\{\phi_n\}$ is  $T$-${\cal G}_n$ full.

Write $T(a)=(T(a)_\N, T(a)_\R)$ for all $a\in A_+\setminus \{0\}.$
Suppose that $a\in A_+\cap (\cup_{n=1}^{\infty}{\cal G}_n)$ is a
non-zero element. It is easy to see that there are
$x_1,x_2,...,x_N\in q_{\infty}(\{B_n\})$ with $\|x_i\|\le
T(a)_\R,$ $i=1,2,...,N$ and $N\le T(a)_\N$ such that
$$
\sum_{i=1}^Nx_i ^*\Phi(a) x_i=1_{q_{\infty}(\{B_n\})}.
$$
Thus $\Phi(a)$ does not belong to any proper ideal of
$q_{\infty}(\{B_n\}).$

Suppose that $b\in A_+\setminus \{0\}$ such that $\Phi(b)\in I$
for some (closed two-sided ideal) I of $q_{\infty}(\{B_n\}).$ Let
$J=\{b\in A: \Phi(b)\in I\}.$ Then $J$ is a (closed two-sided)
ideal of $A$ and for any non-zero $a\in A_+\cap
(\cup_{n=1}^{\infty}{\cal G}_n),$ $a\not\in J.$ Let $b\in
J_+\setminus \{0\}$ with $\|b\|=1.$ Let $0<1/n<1/16.$ There exists
$c\in A_+\cap {\cal G}_n'$ with $\|c\|\ge 1-1/8n$ such that
$$
\|c-b\|<1/2n.
$$
By Lemma 2.7 of \cite{hlx}, there is $x\in A$ such that
$$
x^*x=a=f_{1/n}^{1-1/n}(c) \andeqn xx^*\in \overline{bAb}.
$$
Therefore $x^*x\in J.$ This contradicts with the fact that
$a\not\in J.$

\end{proof}

\begin{prop}\label{fullP1}
Let $A$ be a unital separable amenable \CA\, and let $B$ be a
unital \CA. Suppose that $\{\phi_n\}$ is a sequence of \morp s
from $A$ to $B$ such that
$$
\lim_{n\to\infty}\|\phi_n(a)\phi_n(b)-\phi_n(ab)\|=0\tforal a,
b\in A.
$$
Then there exists a map $T: A_+\setminus \{0\}\to \N\times \R_+$
such that $\{\phi_n\}$ is a $T$-full asymptotic embedding if and
only if $\pi\circ (\{\phi_n\})$ is full.

\end{prop}

\begin{proof}
Suppose that $T: A_+\setminus \{0\}\to \N\times \R_+$ is a map and
$\{\phi_n\}$ is a $T$-full asymptotic embedding. Then, for any
$a\in A_+\setminus \{0\},$ let $N=T(a)_\N$ and $R=T(a)_\R.$ Then
there are $x_{n,i}\in A$ with $\|x_{n,i}\|\le R,$ $i=1,2,...,N$
such that
$$
\sum_{i=1}^N x_{n,i}^*\phi_n(a)x_{n,i}=1_B.
$$
Write $X_i=\{x_{n,i}\}\in l^{\infty}(\{B\})$ (with $\|X_i\|\le
R$). Then
$$
\sum_{i=1}^N \pi(X_i)^*\pi\circ
(\{\phi_n\})(a)\pi(X_i)=1_{q_{\infty}(\{B\})}.
$$
It follows that $\pi\circ (\{\phi_n\})$ is full.

Conversely, if $\pi\circ (\{\phi_n\})$ is full, then for any $a\in
A_+\setminus \{0\},$ there is $N\ge 1$ and $R>0,$ and there are
$y_i\in q_{\infty}(\{B\})$ such that $\|y_i\|\le R,$
$i=1,2,...,N,$ and
$$
\sum_{i=1}^Ny_i^*\pi\circ
(\{\phi_n(a)\})y_i=1_{q_{\infty}(\{B\})}.
$$
There is $X_i=\{x_{n,i}\}\in l^{\infty}(\{B\})$ such that
$\|X_i\|\le R$ and $\pi(X_i)=y_i,$ $i=1,2,....$ Then,
$$
\lim_{n\to\infty}\|\sum_{i=1}^N x_{n,i}^*\phi_n(a)x_{n,i}-1_B\|=0
$$
for all $a\in A_+\setminus \{0\}.$ Thus, for any $a\in
A_+\setminus \{0\},$ there exits an integer $n_a>0$ such that
$$
\|\sum_{i=1}^N x_{n,i}^*\phi_n(a)x_{n,i}-1_B\|<1/2
$$
for all $n\ge n_a.$ One obtains an invertible  positive element
$b_n\in B_+$ with $\|b_n\|\le 2R$ for all $n\ge n_a$ such that
$$
\sum_{i=1}^R b_nx_{n,i}^*\phi_n(a) x_{n,i}b_n=1_B
$$
for all $n\ge n_a.$

Define $T: A_+\setminus\{0\}\to \N\times \R_+$ by $T(a)_\N=N$ and
$T(a)_\R=2R.$ Then $\{\phi_n\}$ is $T$-full.

\end{proof}

\begin{lem}\label{Lunq}
Let $A$ be a separable unital amenable \CA\, in ${\cal N}.$
For any $\ep>0,$ any finite subset ${\cal F}\subset A$ and any $T:
A_+\setminus \{0\}\to \N\times \R_+,$
 there exists a finite subset
${\cal G}={\cal G}(\ep, {\cal F}, A)\subset A,$ $\dt=\dt(\ep,
{\cal F}, A)>0,$ and a finite subset ${\cal P}={\cal P}(\ep, {\cal
F}, A)\subset K_0(A),$ an integer $N=N(\ep, {\cal F}, A)>0$
and ${\cal S}={\cal S}(\ep, {\cal F}, A,
T)\subset A_+$  satisfying the following: for any  unital AF-algebra $F$ with
divisible $K_0(F)$ and any
 $T$-${\cal S}$-full a\morp\, $f: A\to F$ with
\beq\label{Lu-1}
\|f(a)f(b)-f(ab)\|<\dt\andeqn \|f(a)\|\ge 1/2\|a\|
\eneq
for all $a, b\in {\cal G}$  and any two  unital ${\cal
G}$-$\dt$-multiplicative \morp s $\phi, \psi: A\to F$ with
\beq\label{Lu2}
[\phi]|_{\cal P}=[\psi]|_{\cal P}
\eneq
there exists  a  unitaries $u\in M_{N+1}(F)$ such that,
\beq\label{Lu-3}
\|{\rm ad}\, u\circ {\rm diag}(\phi(a),
\overbrace{f(a),f(a),...,f(a)}^N)- {\rm
diag}(\psi(a),\overbrace{f(a),f(a),...,f(a)}^N)\|<\ep
\eneq
for all $a\in {\cal F}.$

Moreover, the assumption that $K_0(F)$ is divisible can be
replaced by $K_1(A)$ is torsion free and
\beq\label{Lu2+}
[f]|_{\cal P}=[\phi]|_{\cal P}.
\eneq

\end{lem}

\begin{proof}
Suppose that the lemma is false. Then there are $\ep_0>0$ and a
finite subset ${\cal F}_0\subset A$ and $T: A_+\setminus\{0\} \to
(\N, \R_+)$  such that there are a sequence of decreasing positive
numbers $\dt_n\to 0,$ an increasing sequence of finite subsets
$\{{\cal G}_n\}\subset A$ whose   union is dense in $A,$ an
increasing sequence of finite subsets $\{{\cal P}_n\}\subset
K_0(A)$ whose union is $K_0(A),$ an increasing sequence $\{k(n)\}$
of positive integers for which $\lim_{n\to\infty} k(n)=\infty,$ an
increasing sequence of finite subsets $\{S_n\}$ which satisfies
the conclusion of \ref{Ldense}, a sequence $\{F_n\}$ of  unital AF
\CA s with divisible $K_0(F_n),$ and three sequences of ${\cal
G}_n$-$\dt_n$-multiplicative \morp s $\{\phi_n\},$ $\{\psi_n\}$
and $\{f_n\}$ from $A$ into $F_n$ for which $\{f_n\}$ is
$T$-${\cal S}_n$ full and for which
\beq\label{Lu-4-1}
[\phi_n]|_{{\cal P}_n}=[\psi_n]|_{{\cal P}_n},
\eneq
satisfying:
\beq\label{1nL1}
\hspace{-0.2in}\inf\{\sup\{\|u^*{\rm diag}(\phi_n(a),\overbrace{
f_n(a),...,f_n(a)}^{k(n)})u -{\rm
diag}(\psi_n(a),\overbrace{f_n(a),...,f_n(a)}^{k(n)})\|: a\in
{\cal F}_0\}\ge \ep_0,
\eneq
where the infimum is taken over all unitaries in
$M_{k(n)+1}(B_n).$

 Let $\pi: l^{\infty}(\{F_n\})\to q_{\infty}(\{F_n\})$ be
the quotient map. Put ${\bar \Phi}=\pi\circ \{\phi_n\},$ ${\bar
\Psi}=\pi\circ\{\psi_n\}$ and ${\bar H}=\pi\circ \{f_n\}.$  Then
${\bar \Phi},$ ${\bar \Psi}$ and ${\bar H}$ are \hm s. By
\ref{Ldense}, ${\bar H}$ is a full monomorphism from $A$ to
$q_{\infty}(\{F_n\}).$ Moreover, by (\ref{Lu-4-1}), since each
$F_n$ is AF (see Cor. 2.1  of \cite{GL}),
\beq\label{Lu-4-2}
{\bar \Phi}_{*0}={\bar \Psi}_{*0}.
\eneq

Since $F_n$ is AF, $K_1(F_n)=\{0\}$ and $K_0(F_n)$ is torsion
free. By Cor. 2.1 of \cite{GL}, since $F_n$ are AF,
\beq\label{Lu-4-K5}
K_1(q_{\infty}(\{F_n\}), \Z/k\Z)=\{0\},\,\,\,k=2,3,....
\eneq
Moreover (by Cor. 2.1 of \cite{GL}),
\beq\label{Lu-4-K6}
K_1(q_{\infty}(\{F_n\}))=\{0\}.
\eneq

Since $K_0(F_n)$ is divisible,
\beq\label{Lu-4-K1}
\prod_n K_0(F_n,\Z/k\Z)/\oplus_n K_0(F_n,
\Z/k\Z)=\{0\},\,\,\,i=0,1,\,k=2,3,....
\eneq
By Cor. 2.1 of \cite{GL},
\beq\label{Lu-4-K2}
K_0(q_{\infty}(\{F_n\}),\Z/k\Z)=\{0\}, \,i=0,1,\, k=2,3,...
\eneq
Since $A$ is assumed to satisfy the UCT, we compute that
\beq\label{Lu-4-K3}
[{\bar \Phi}]=[{\bar \Psi}]\,\,\,{\rm in}\,\,\,KL(A,
q_{\infty}(\{F_n\})).
\eneq

If the condition  that $K_0(F_n)$ is divisible is replaced by
$K_1(A)$ is torsion free and $[f_n]|_{{\cal P}_n}=[\phi_n]|_{{\cal
P}_n},$ then we should get
\beq\label{Lu-4-K4}
{\bar H}_{*0}={\bar \Phi}_{*0}={\bar \Psi}_{*0} \andeqn K_0(A,
\Z/k\Z)=K_0(A)/kK_0(A),\,\,k=2,3,...
\eneq

Let ${\cal Q}\subset \underline{K}(A)$ be a finite subset. There
is an integer $k>0$ such that
\beq\label{Lu-6}
{\cal Q}\cap K_i(A,\Z/m \Z)=\{0\}\rforal m\ge k,
\eneq
$i=0,1.$ Define ${\tilde H}=\overbrace{{\bar H}\oplus {\bar
H}\oplus \cdots \oplus {\bar H}}^{k!-1}.$ By (\ref{Lu-4-K4}),
(\ref{Lu-4-K5}) and (\ref{Lu-4-K6}), we compute that
\beq\label{Lu-7}
[{\bar \Phi}\oplus {\tilde H}]|_{{\cal Q}\cap K_i(A,
\Z/m\Z)}=[{\bar \Psi}\oplus {\tilde H}]|_{{\cal Q}\cap
K_i(A,\Z/m\Z)},
\eneq
if $i=2,3,...,k.$
%
By (\ref{Lu-4-K6}),
we also have that
\beq\label{Lu-8}
[{\bar \Phi}\oplus {\tilde H}]|_{{\cal Q}\cap K_1(A)}=[{\bar
\Psi}\oplus {\tilde H}]|_{{\cal Q}\cap K_1(A)}.
\eneq
We conclude that
\beq\label{Lu-9}
[{\bar \Phi}\oplus {\tilde H}]]|_{\cal Q}=[{\bar \Psi}\oplus
{\tilde H}]]|_{\cal Q}.
\eneq
Using either (\ref{Lu-4-K3}) or (\ref{Lu-9}) with sufficiently
large ${\cal Q}$ (which depends on $\ep_0$ and ${\cal F}_0$) (and
hence large $k$), and by applying  3.9 of \cite{Lnuct}  we obtain
an integer $K>0$ and a unitary $U\in
M_{2K+1}(M_{k!}(q_{\infty}(\{F_n\})))$ such that
\beq\label{Lu-10}
\hspace{-0.2in}\|U^*{\rm diag}({\bar \Phi}(a)\oplus {\tilde
H}(a),\overbrace{{\tilde H}(a),\cdots {\tilde H}(a)}^{2K})U -{\rm
diag}({\bar \Psi}(a)\oplus {\tilde H}(a),\overbrace{{\tilde
H}(a),\cdots {\tilde H}(a)}^{2K})\|<{\ep_0\over{4}}
\eneq
for all $a\in {\cal F}.$ Choose $N=(2K+1)(k!-1).$ It is easy to
see that there is a unitary $W=\{w_n\}\in l^{\infty}(\{B_n\})$
such that $\pi(W)=U.$ Therefore, there is $N_0>0$ such that, for
any $n\ge N_0,$
\beq\label{1nL3}
\|u_n^*{\rm diag}(\phi_n(a),\overbrace{f_n(a),...,f_n(a)}^{N})u-
{\rm diag}(\psi_n(a),
\overbrace{f_n(a),...,f_n(a)}^{N})\|<{\ep_0\over{2}}
\eneq
for all $a\in {\cal F}_0.$

The above estimate (\ref{1nL3}) contradicts with (\ref{1nL1}).

\end{proof}

\begin{df}\label{dCH}
{\rm Denote by ${\cal C}_0$ the class of unital \CA\, $C$  with
the form $C=M_n(C(X)),$ where $X$ is a compact subset of a  finite CW complex,
$n$ is a positive integer, $P\in
M_n(C(X))$ is a projection. Note that $X$ is not assumed to be connected.

Recall that a \CA\, is said to be AH, if it is an inductive limit
of finite direct sums of \CA s in ${\cal C}_0.$

}
\end{df}

The following follows from the results in \cite{Lncd} immediately.

\begin{thm}\label{cdT1}
Let $C$ be a unital AH-algebra and let $A$ be a unital simple
\CA\, with tracial rank zero. Suppose that $h_1, h_2: C\to A$ are
two monomorphisms such that
\beq\label{cdt1}
[h_1]=[h_2]\,\,\,\text{in}\,\,\, KL(C,A)\andeqn
\eneq
\beq\label{cdt2}
\tau\circ h_1(f)=\tau\circ h_2(f)\tforal f\in C\andeqn \tau\in
T(A).
\eneq
Then $h_1$ and $h_2$ are approximately unitarily equivalent,
i.e., there exists a sequence of unitaries $\{u_n\}\subset A$
such that
\beq\label{cdt3}
\lim_{n\to\infty}{\rm ad}\, u_n\circ h_1(f)=h_2(f)\tforal f\in C.
\eneq
\end{thm}

\begin{proof}
First, we note that it follows from Cor. 4.8   of \cite{Lncd}  (see also 3.4 of \cite{Lncd}) that the theorem holds for $C\in {\cal C}_0.$
 The case that
$C$ is a finite direct sum of \CA s in ${\cal C}_0$ then also
follows.

 Now we assume that $C=\lim_{n\to\infty} (A_n, \phi_n),$ where
 each $A_n$ is a finite direct sum of \CA s in ${\cal C}_0.$
Since the quotients of \CA s in ${\cal C}_0$ are still in ${\cal
C}_0,$ without loss of generality, we may assume that each
$\phi_n$ is injective. Denote by $\phi_{n, \infty}: A_n\to C$ the
embedding induced by the inductive limit (note we now assume that
each $\phi_n$ is injective).

Let $\ep>0$ and let ${\cal F}\subset C$ be a finite subset.
Without loss of generality, we may assume that ${\cal F}\subset
\phi_{n, \infty}(A_n).$ Denote by ${\cal F}'\subset A_n$ the
finite subset so that $\phi_{n,\infty}({\cal F}')={\cal F}.$

Now
$$
[h_1\circ \phi_{n, \infty}]=[h_2\circ \phi_{n, \infty}]\,\,\, \text{in}\,\,\,
KL(A_n, A) \andeqn
$$
$$
\tau\circ h_1\circ \phi_{n, \infty}(f)=\tau\circ h_2\circ \phi_{n,
\infty}(f)
$$
for all $f\in A_n$ and $\tau\in T(A).$ By what we have shown (by
applying \cite{Lncd}) for the cases that \CA s are direct sum of
\CA s in ${\cal C}_0,$ we obtain a unitary in $A$ such that
$$
{\rm ad}\, u\circ h_1\circ \phi_{n, \infty}\approx_{\ep} h_2\circ
\phi_{n,\infty}\,\,\,\text{on}\,\,\, {\cal F'}.
$$
It follows that
$$
{\rm ad}\, u\circ h_1\approx_{\ep}h_2\,\,\,\text{on}\,\,\, {\cal
F}.
$$

\end{proof}

\section{AH-algebras}

\begin{thm}\label{embAH}
Let $A$ be a unital  AH-algebra. Then $A$ can be embedded into a unital
simple AF-algebra if and only if $A$ has a faithful tracial state.

Moreover, if $t$ is a faithful tracial state, then there is a
unital simple AF-algebra $B$ with a unique tracial state $\tau$
such that $(K_0(B), K_0(B)_+, [1_B])=(\D, \D_+, 1),$ where $\D$ is
a divisible countable dense additive subgroup of $\R,$
$K_0(B)=\rho_B(K_0(B))$ and there is a unital monomorphism $\phi:
A\to B$ such that
$$
\tau\circ \phi=t.
$$
\end{thm}

\begin{proof}
Let $t$ be a faithful tracial state on $A.$ Write
$A=\lim_{n\to\infty} (A_n, h_n),$ where\\
$A_n=\oplus_{i=1}^{j(n)} P_{n,j} M_{d(i,n)}(C(X_{n,i}))P_{n,j},$
where each $d(i,n)$ is a positive integer, $X_{n,i}$ is a compact subset of a finite CW complex
and $P_{n,j}\in M_{d(i,n)}(C(X_{n,i})$ is a
projection and each $h_n$ is injective. Moreover, we assume that
there are integers $r(n,i), s(n,i)>0$ and a projection
$$
Q_{n,i}\in M_{s(n,i)}(P_{n,j} M_{d(i,n)}(C(X_{n,i}))P_{n,j})
$$
such that
$$
Q_{n,i}M_{s(n,i)}(P_{n,j}
M_{d(i,n)}(C(X_{n,i}))P_{n,j})Q_{n,i}\cong M_{r(n,i)}(C(X_{n,i})).
$$
We extend the monomorphism $h_{n,\infty}: M_{s(n,i)}(P_{n,j}
M_{d(i,n)}(C(X_{n,i}))P_{n,j}) \to M_{s(n,i)}(A).$ To simplify
notation, we may identify $Q_{n,i}M_{s(n,i)}(P_{n,j}
M_{d(i,n)}(C(X_{n,i}))P_{n,j})Q_{n,i}$ with
$M_{r(n,i)}(C(X_{n,i})).$ Let $e_{11}^{(n,i)}$ be a rank one
projection in $M_{r(n,i)}(C(X_{n,i})).$ Consider the monomorphism
$h_{n,\infty}^{(i)}:
e_{11}^{(n,i)}M_{r(n,i)}(C(X_{n,i}))e_{11}^{(n,i)} \to
h_{n,\infty}(e_{11}^{(n,i)})M_{s(n,i)}(A)h_{n,\infty}(e_{11}^{(n,i)})$
which is the restriction of $h_{n, \infty}.$
Note that
$$C(X_{n,i})\cong e_{11}^{(n,i)}M_{r(n,i)}(C(X_{n,i}))e_{11}^{(n,i)}.$$
 We
use $t$ for $t\otimes Tr,$ where $Tr$ is the standard trace on
$M_{r(n,i)}.$ This gives a strictly positive probability Borel
measure $\mu_{n,i}$ on $X_{n,i}.$ As in  3.2 of \cite{Lnemb2}, one
obtains a sequence of partitions ${\cal P}_m:$ $Y_{n,i}(m,j),$
$j=1,2,...c(m),$ which has the following properties:

(1) $Y_{n,i}(m,j)\cap Y_{n,i}(m,j')=\emptyset,$ if $j\not=j',$

(2) $\lim_{m\to\infty}\max\{{\rm diam}(Y_{n,i}(m,j)):
j=1,2,...,c(m)\}=0,$

 (3) ${\rm int}(Y_{n,i}(m,j))\subset Y_{n,i}(m,j)\subset
\overline{Y_{n,i}(m,j)}$ and each $Y_{n,i}(m,j)$ is a Borel set,

(4) $\mu_{n,i}(\partial(Y_{n,i}(m,j))=0,$

(5) $X_{n,i}=\cup_{j=1}^{c(m)}Y_{n,i}(m,j).$

Moreover, $Y_{n,i}(m,j)$ is a union of finitely many members of
$Y_{n,i}(m+1,j').$

Put $D_{n,i}=\{\mu_{n,i}(Y_{n,i}(m,j)): m,j\}.$ Let $\D$ be the
additive subgroup of $\R$ generated by $\Q$ and all $D_{n,i}.$
Then $\D$ is a divisible countable dense subgroup of $\R.$ Let $C$
be a unital simple AF-algebra with $(K_0(C), K_0(C), [1_C])=(\D,
\D_+, 1).$ Let $\tau$ be the unique tracial state on $C.$

Fix $n,$ let $R(n)=\max\{r(n,i): i=1,2,...,d(n)\}.$ There are
mutually orthogonal projections $p_{n,i}\in M_{R(n)}(C)$ such that
\beq\label{AH1}
\tau(p_{n,i})=t\circ
h_{n,\infty}(e_{11}^{(n,i)}),\,\,\,i=1,2,...,d(n).
\eneq
Fix $i$ and $m,$  there are mutually orthogonal projections
$\{q_{n,i}(m,j)\} \subset p_{n,i}Cp_{n,i}$ such that
\beq\label{AH2}
\sigma(e_{n,i}(m,j))=\mu_{n,i}(Y_{n,i}(m,j)),\,\,\, j=1,2,...,c(m).
\eneq
If, for some finite subset $J,$ $\sum_{j'\in J} Y_{n,i}(m+1,
j')=Y_{n,i}(m,j),$ then there are mutually orthogonal projections
$e_{n,i}(m+1,j')\in e_{n,i}(m,j)Ce_{n,i}(m,j)$ such that
\beq\label{AH2+1}
\sum_{j'\in J} e_{n,i}(m+1, j')&=&e_{n,i}(m,j)\andeqn\\
\tau(e_{n,i}(m+1,j'))&=&\mu_{n,i}(Y_{n,i}(m+1,j')).
\eneq

Therefore
$$
\{e_{n,i}(m,j): j=1,2,...,c(m), m=1,2,...\}
$$
is a set of mutually orthogonal projections. Let $C(n,i)\subset
p_{n,i}Cp_{n,i}$ be the \SCA\, generated by $\{e_{n,i}(m,j): 1\le
j\le c(m), m\in \N\}.$

Define $h_m: C(X_{n,i})\to p_{n,i}Cp_{n,i}$ by
\beq\label{AH3}
h_m^{(n,i)}(f)=\sum_{j=1}^{c(m)} f(\xi_{n,i,m,j})e_{n,i}(m,j)
\eneq
for $f\in C(X_{n,i}),$ where $\xi_{n,i, m,j}\in Y_{n,i}(m,j)$ is a
point. By (2) above and (\ref{AH2+1}), for each $n$ and $i,$ $\{h_m^{(n,i)}(f)\}$ is
a Cauchy sequence. Define
\beq\label{AH4}
h^{(n,i)}(f)=\lim_{m\to\infty} h_m^{(n,i)}(f)
\eneq
for $f\in C(X_{n,i}).$ We extend, for each $n$ and $i,$
$h^{(n,i)}: M_{r(n,i)}(C(X_{n,i)}))\to
M_{r(n,i)}(p_{n,i}Cp_{n,i}).$
By (\ref{AH1}), (\ref{AH2}) and (\ref{AH3}) as well as (\ref{AH4}),
 we obtain a monomorphism $h^{(n)}:
A_n\to M_{R(n)}(C)$ such that
\beq\label{AH5}
\tau\circ h^{(n)}(a)=t\circ h_{n, \infty}(a)
\eneq
for all $a\in A_n.$ Since
\beq\label{AH6}
\tau(h^{(n)}(1_{A_n}))\le 1\andeqn {\rm ker}\rho_C=\{0\},
\eneq
there is a unitary $v_n\in M_{R(n)}(C)$ such that ${\rm
ad}\,v_n\circ h^{(n)}$ maps $A_n$ into $C$ such that
\beq\label{AH7}
\tau\circ {\rm ad}\, v_n\circ h^{(n)}(a)=t\circ h_{n,\infty}(a)
\eneq
 for all $a\in A_n.$
 It follows that, for any $n'>n,$
\beq\label{AH8}
\tau\circ {\rm ad}\,v_{n'}\circ h^{(n')}\circ h_{n,n'}(a)= t\circ
h_{n', \infty}\circ h_{n,n'}(a).
\eneq
This implies (${\rm ker}\rho_C=\{0\}$) that
\beq\label{AH9}
({\rm ad}\,v_{n'}\circ h^{(n')}\circ h_{n,n'})_{*0}=({\rm ad}\,
v_n\circ h^{(n)})_{*0}.
\eneq
Since $K_1(C)=\{0\},$ $K_0(C)$ is torsion free and divisible, we
have
\beq\label{AH10}
[{\rm ad}\,v_{n'}\circ h^{(n')}\circ h_{n,n'}]=[{\rm ad}\,
v_n\circ h^{(n)}]\,\,\,{\rm in}\,\,\, KL(A_n,C).
\eneq


Combining (\ref{AH10}) and (\ref{AH8}), by applying \ref{cdT1}
and passing to a subsequence, if necessary, we obtain unitaries
$u_n\in C$ such that
\beq\label{AH11}
{\rm ad}\, u_{n+1}\circ {\rm ad}\, v_{n+1}\circ h^{(n+1)}\circ
h_{n}\approx_{1/2^n} {\rm ad}\, v_n\circ h^{(n)}\,\,\,{\rm
on}\,\,\, {\cal F}_n.
\eneq
Put $\phi_n={\rm ad}\, u_n\circ {\rm ad}\, v_n\circ h^{(n)},$
$n=1,2,....$ Thus we obtain the following
 approximate intertwining:
$$
\begin{array}{ccccccc}
A_1 & {\stackrel{h_1}{\to}} & A_2  &{\stackrel{h_2}{\to}} & A_3 & {\stackrel{h_3}{\to}} \cdots & A\\
\downarrow_{\phi_1} && \downarrow_{\phi_{2}} & &\downarrow_{\phi_3} & \cdots \\
\hspace{-0.1in}C & {\stackrel{{\rm id}_C}{\to}} &
\hspace{-0.1in}C & {\stackrel{{\rm id}_C}{\to}} & \hspace{-0.1in}C
& {\stackrel{{\rm id}_C}{\to}} \cdots  &C
\end{array}
$$

Thus, since each $\phi_n$ is injective,  we obtain a unital embedding $\phi: A\to C$ such that
$\tau\circ \phi(a)=t(a)$ for all $a\in A.$

\end{proof}

\begin{cor}\label{CembAH}
Let $A$ be a unital AH-algebra and let $\af\in Aut(A).$ Suppose
that $A$ admits a faithful $\af$-invariant tracial state $t.$ Then
there is a unital simple AF-algebra $B$ with a unique tracial
state $\tau$ and with $(K_0(B), K_0(B)_+, [1_B])=(\D, \D_+, 1),$
where $\D=\rho_B(K_0(B)),$ and there is a unital embedding $\phi:
A\to B$ such that
$$
\tau\circ \phi\circ \af=\tau\circ \phi.
$$
Furthermore, we may assume that $\D$ is divisible.
\end{cor}

\begin{proof}
In the proof of \ref{embAH}, we may begin with the tracial state
$t.$ Thus
$$
\tau\circ \phi\circ \af=t\circ \af=t=\tau\circ \phi.
$$

\end{proof}

The following is certainly known.

\begin{cor}\label{AFemb}
Let $A$ be a unital AF-algebra. Then $A$ can be embedded into a
unital simple AF-algebra if and only if $A$ admits a faithful
tracial state.
\end{cor}

\section{Crossed products of AH-algebras by $\Z$}

The following is proved in \cite{Lnhomp}. The constant $\dt$ and
the finite subset ${\cal G}$ can be chosen independent of $A$ and
$h$ but  they are dependent on a measure distribution. Please see
2.10 of \cite{Lnhomp} for the definition of
$\text{Bott}(h,v)|_{\cal P}.$  See 2.10  and 2.11 of \cite{Lnhomp}
also for the definition of $\text{bott}_1(h,v)$ which is  related
to the usual bott element $\text{bott}_1(u,v)$ for a pair of
unitaries.

\begin{thm}{\rm (\cite{Lnhomp})}\label{BHPT}
Let $A$ be a unital AH-algebra and let $B$ be a unital separable
simple \CA\, of tracial rank zero. Suppose that $h: A\to B$ is a
unital monomorphism. Then, for any $\ep>0$ and any finite subset
${\cal F}\subset A,$ there exists a positive number $\dt,$ a
finite subset ${\cal G}\subset A$ and a finite subset ${\cal
P}\subset \underline{K}(A)$ satisfying the following: suppose that
$w, u\in U(B)$ are unitaries  such that
$$
\|{\rm ad}\, w\circ h(a)u-u{\rm ad}\, w\circ h(a)\|<\dt\tforal
a\in {\cal G}\,\,\,and\,\,\, \rm{Bott}({\rm ad}\,w\circ
h,u)|_{\cal P}=0.
$$
Then there is a continuous rectifiable path $\{u_t: t\in [0,1]\}\subset U(B)$
such that
$$
u_0=u,\,\,\,u_1=1_B\,\,\, and
$$
$$
\|{\rm ad}\, w\circ h(a)u_t-u_t{\rm ad}\, w\circ h(a)\|<\ep\tforal
a\in {\cal F}\,\,\, and  \tforal t\in [0,1].
$$
Moreover
$$
\|u_t-u_{t'}\|\le (2\pi+\ep)|t-t'|\tforal t,t'\in [0,1].
$$
In particular,
$$
\text{Length}(\{u_t\})\le 2\pi+\ep.
$$
\end{thm}

Using the above theorem, we obtain the following:

\begin{lem}\label{path}
Let $A$ be a unital AH-algebra and let $B$ be a unital separable
simple \CA\, of tracial rank zero. Suppose that $h: A\to B$ is a
unital monomorphism.

Then, for any $\ep>0$ and any finite subset ${\cal F}\subset A,$
there exists $\dt>0$ and a finite subset ${\cal G}\subset A$
satisfying the following: If there is a continuous path of
unitaries $\{u_t: t\in [0,1]\}$ such that $u_1=1_B,$
$$
\|h(b)u_t-u_th(b)\|<\dt\tforal b\in {\cal G}\,\,\,and\,\,\, t\in [0,1],
$$
then there is a continuous path of unitaries $\{w_t: t\in [0,1]\}$
such that
$$
w_0=u_0,\,\,\,w_1=1_B,
$$
$$
\|h(a)w_t-w_t h(a)\|<\ep\tforal a\in {\cal F}\,\,\,and \,\,\, \tforal t\in
[0,1].
$$
Moreover,
$$
\|u_t-u_{t'}\|\le (2\pi+\ep)|t-t'|\tforal t,t'\in [0,1].
$$
\end{lem}

\begin{proof}
Fix $\ep$ and a finite subset ${\cal F}\subset A.$ Without loss of
generality, we may assume that ${\cal F}$ is in the unit ball of
$A.$ Let $\dt_1>0$ and ${\cal G}_1\subset A$ and ${\cal P}\subset
\underline{K}(A)$ be finite subsets required in \ref{BHPT} for the
above $\ep$ and ${\cal F}.$ We choose $\dt_2>0$ and a finite
subset ${\cal G}_2'\subset A\otimes C(S^1)$ so that for any
$\dt_2$-${\cal G}_2'$-multiplicative \morp\, $L_1, L_2: A\otimes
C(S^1)\to B,$ $[L_i]|_{{\boldsymbol{\bt}}(\cal P)}$ is well defined for
$i=1,2$ and
$$
[L_1]|_{{\boldsymbol{\bt}}(\cal P)}=[L_2]|_{{\boldsymbol{\bt}}(\cal P)},
$$
provided that
$$
L_1\approx_{\dt_2} L_2\,\,\,\text{on}\,\,\,{\cal G}_2'.
$$
(see 2.4 of \cite{Lnhomp} and 2.10 of \cite{Lnhomp} for the
definition of Bott map ${\boldsymbol{\bt}}$).

Without loss of generality, by choosing smaller $\dt_2,$ we may
assume that ${\cal G}_2'={\cal G}_2\otimes S,$ where
$S=\{1_{C(S^1)}, z\}$ and $z$ is the identity function on the unit
circle.

We may assume that ${\cal F}\cup {\cal G}_1\subset {\cal G}_2$ and
$\dt_2<\min\{\ep, \dt_1\}.$

Let $\dt_3>0$ (in place of $\ep$) and ${\cal G}_3\subset A$ (in
place of ${\cal F}_1$) be a finite subset required by 2.7 of
\cite{Lnhomp} for $\dt_2/2$ (in place of $\ep$) and ${\cal G}_2$
(in place of ${\cal F}_2$).

Let $\dt=\min\{\dt_3, \dt_2\}$ and ${\cal G}={\cal G}_3\cup {\cal
G}_2.$

We may assume that ${\cal G}_3$ is a subset of the unit ball of
$A.$

Now suppose that $\{w_t: t\in [0,1]\}$ and $h$ are in the lemma
for the above $\dt$ and ${\cal G}.$

It follows from 2.7 of \cite{Lnhomp} that there is, for each $t,$
an $\dt_2$-${\cal G}_2$-multiplicative \morp\, $\psi_t: A\otimes
C(S^1)\to B$ such that
\beq\label{path-1}
\|\psi_t(a\otimes 1)-h(a)\|<\dt_2/2\andeqn \|\psi_t(a\otimes g)-
h(a)g(u_t)\|<\dt_2/2
\eneq
for all $a\in {\cal G}_2$ and $g\in S.$

There is a partition $0=t_0<t_1<\cdots <t_m=1$ such that
\beq\label{path-2}
\|u_{t_{i-1}}-u_{t_i}\|<\dt_2/2,\,\,\,i=1,2,...,m.
\eneq
 It follows from (\ref{path-1}) and (\ref{path-2}) that
 \beq\label{path-3}
 \|\psi_{t_{i-1}}(a\otimes g)-\psi_{t_i}(a\otimes
 g)\|<\dt_2\rforal a\in {\cal G}_2\andeqn g\in S.
 \eneq
 Moreover,
\beq\label{path-4}
\|\psi_{t_m}(a\otimes g)-h(a)\|<\dt_2\rforal a\in {\cal
G}_2\andeqn g\in S.
\eneq
It follows (from the choice of $\dt_2$ and ${\cal G}_2$) that
\beq\label{path-5}
[\psi_{t_m}]_{{\boldsymbol{\bt}}({\cal P})}=0.
\eneq

 Therefore, by the choice of $\dt_2$ and ${\cal G}_2,$ and by
 (\ref{path-5}),
\beq\label{path-6}
[\psi_{t_0}]|_{{\boldsymbol{\bt}}({\cal P})}=0.
\eneq
In other words,
\beq\label{path-7}
\text{Bott}(h, u_0)|_{\cal P}=0.
\eneq
Since $\dt<\dt_1$ and ${\cal G}_1\subset {\cal G}_2,$ by applying
\ref{BHPT}, we obtain a continuous path of unitaries $\{w_t: t\in
[0,1]\}$ such that
$$
u_0=u_0,\,\,\, u_1=1_B\andeqn
 \|h(a)u_t-u_th(a)\|<\ep\rforal a\in {\cal F}\andeqn t\in [0,1].
 $$
 Moreover,
 $$
 \|u_t-u_{t'}\|\le (2\pi+\ep)|t-t'|\tforal t, t'\in [0,1].
 $$

\end{proof}

\begin{df}\label{LlL}
Let $A$ be a unital separable \CA\, and let $B$ be a unital \CA.
Suppose that $h: A\to B$ is a unital monomorphism.
We say $h$ satisfies property (H) for a positive number $L,$  if the following holds:
For any $\ep>0$ and any finite subset ${\cal F}\subset A,$ if there is a continuous path of
unitaries $\{v_t: t\in [0,1]\}$ in $B$ with $v(1)=1$ such that
$$
\|v_t {\rm ad}\, w\circ h(a)-{\rm ad}\,w\circ h(a)v_t\|<\dt\tforal a\in {\cal G}\andeqn t\in [0,1],
$$
where $w\in B$ is a unitary,
there is a continuous path of unitaries $\{u_t: t\in [0,1]\}$ such that
$$
u_0=v_0,\,\,\, v_1=1\,\,\,\andeqn \|u_t{\rm ad}\, w\circ h(a)-{\rm ad}\,w\circ h(a)u_t\|<\ep
$$
for all $a\in {\cal F}$ and all $t\in [0,1].$ Moreover,
$$
\|u_t-u_{t'}\|\le L|t-t'|\tforal t, t'\in [0,1].
$$
\end{df}

It should be noted that a rectifiable continuous path with length, say $2\pi+\ep,$ may not
have the Lipschitz condition.

The following lemma is taken from \cite{M1}. The proof is also the same as
in \cite{M1}.
However, since it has  different assumption and slightly different conclusion, for the convenience,
we will present the  proof here.

\begin{lem}{\rm (Proposition 1 of \cite{M1})}\label{LBM}
Let $A$ be a unital separable \CA, let $B$ be a unital simple
separable \CA\, of tracial rank zero, let $\phi: A\to B$ be a
unital embedding, let $\af\in Aut(A)$ and let $\bt_0\in Aut(B)$ be
automorphisms.

Suppose that $\phi\circ \af$ has the property (H) for a positive number $L>0$ and suppose that there is a
continuous path $\{v(t): t\in [0, \infty)\}$  of unitaries in $B\otimes {\cal U}$ satisfying the following:
$$
\lim_{t\to\infty}\|\phi\circ \af(a)-{\rm ad}\, v(t)\circ \bt\circ
\phi(a)\|=0
$$ for all $a\in A,$ where $\bt=\bt_0\otimes \sigma$ and ${\cal U}$ and $\sigma$ are defined in
\ref{DU}.

Then there are  unitaries $w,\, V_n\in B\otimes {\cal U}$ {\rm
(}$n=1,2,...,${\rm )} such that
$$
\phi'\circ\af={\rm ad}\, w\circ \bt\circ  \phi',
$$
where $\phi'(a)=\lim_{n\to\infty}{\rm ad}\, (V_1V_2\cdots
V_n)\circ \phi.$
\end{lem}

\begin{proof}
To simplify notation, we may assume that $B=B\otimes {\cal U}.$
Note that $\bt$ is an approximately inner automorphism on $B$
which satisfies the cyclic Rokhlin property of \ref{DU}.

Let $\{{\cal F}_n\}$  be an increasing sequence of finite subsets
of the unit ball of $A$ such that
$\overline{\cup_{n=1}^{\infty}{\cal F}_n}=A.$ Let $N_1>0$ such
that
\beq\label{LBM-1}
{L+1\over{N_1-1}}<1/16.
\eneq

Since $\phi\circ \af$ has the property (H) for the constant $L>0$ and
for $\ep=1/16$ and ${\cal F}=\cup_{i=0}^{N_1}\af^{-i}({\cal F}_1),$ we let $\eta_1>0$
(in place of $\dt$)
and ${\cal G}_1\subset A$  (in place of ${\cal G}$) be  a finite subset which is required
in \ref{LlL}. We may assume that ${\cal F}_1\subset {\cal G}_1.$

Choose $0<\dt_1<\min\{{\eta_1\over{N_1^2}}, {1\over{32(N_1^2+N_1+1)}}\}.$

By replacing $\bt$ by ${\rm ad} \circ v(t_0)\circ \bt,$ we may assume that
$v(0)=1$ and
\beq\label{LBM-2}
\|\phi\circ \af(a)-{\rm ad}\, v(t)\circ \bt\circ \phi(a)\|<{\dt_1\over{2}}
\eneq
for all $a\in \cup_{i=0}^{N_1-1}\af^{-i}({\cal G}_1)$ and $t\in
[0,\infty).$ In particular,
\beq\label{LBM-2-1}
\|\phi\circ \af(a)-\bt\circ \phi(a)\|<\dt_1/2\,\,\andeqn
\\\label{LBM-2-2}
\|v(t)\phi\circ \af(a)v(t)^*-\phi\circ\af(a)\|<\dt_1/2
\eneq
for all $a\in \cup_{i=0}^{N_1-1}\af^{-i}({\cal G}_1).$

 Choose $N_2>N_1$ such that
${L+1\over{N_2-1}}<{1\over{4\cdot\,4^2}}.$

Let $1>\eta_2>0,$ ${\cal G}_2\subset A$ be a finite subset corresponding to
$1/4^3$ and $\cup_{i=0}^{N_2-1}\af^{-i}({\cal F}_2)$
in \ref{LlL} (for $\phi\circ \af$ and $L>0$).

Choose $\dt_2>0$ such that $\dt_2<\min\{{\eta_2\over{4}}, {1\over{4^2(N_2^2+N_2)}}\}.$

There is  $t_1>0$ such that
\beq\label{LBM-3}
\|\phi\circ \af(a)-{\rm ad}\, v(t)\circ \bt\circ \phi(a)\|<{\dt_2\over{2}}
\eneq
for all $t>t_1$ and $a\in \cup_{i=1}^{N_2-1}\af^{-i}({\cal G}_2).$

For $t\in [0, t_1],$ we define
\beq\label{LBM-4}
u_0(t)&=&1,\,\,\,u_1(t)=v(t),\,\,\, u_2(t)=v(t)\bt(v(t))\\
u_3(t)&=& v(t)\bt(v(t))\bt^2(v(t)),....
\eneq
By ( \ref{LBM-2-1} ),
\beq\label{LBM-6}
\phi\circ \af\approx_{k\dt_1/2}\bt^k\circ \phi\circ \af^{-k+1}\,\,\,\text{on}\,\,\,{\cal G}_1
\eneq
$k=1,2,...,N_1.$
Thus,
\beq\label{LBM-7}
&&\hspace{-0.8in}\|[u_k(t),\,\phi\circ \af(a)]\|= \|[v(t)\bt(v(t))\cdots
\bt^{k-1}(v(t)),\,\phi\circ \af(a)]\|\\
&&\hspace{-0.5in}=\|v(t)\bt(v(t))\cdots \bt^{k-1}(v(t))\phi\circ \af(a)-
v(t)\phi\circ \af(a)\bt(v(t))\cdots \bt^{k-1}(v(t))\|\\
&+&\|v(t)\phi\circ \af(a)\bt(v(t))\cdots \bt^{k-1}(v(t))- \phi\circ
\af(a)v(t)\bt(v(t))\cdots
\bt^{k-1}(v(t))\|\\
&&\hspace{-0.5in}=\|\bt(v(t))\cdots \bt^{k-1}(v(t))\phi\circ \af(a)-
\phi\circ \af(a)\bt(v(t))\cdots \bt^{k-1}(v(t))\|\\
&+&\|v(t)\phi\circ \af(a)- \phi\circ
\af(a)v(t)\|\\
&&\hspace{-0.5in}=\|\bt(v(t))\cdots \bt^{k-1}(v(t))\phi\circ \af(a)-
\bt(v(t))\phi\circ \af(a)\bt^2(v(t))\cdots\bt^{k-1}(v(t))\|\\
&+&\|\bt(v(t))\phi\circ \af(a)\bt^2(v(t))\cdots\bt^{k-1}(v(t))
-\phi\circ \af(a)\bt(v(t))\cdots \bt^{k-1}(v(t))\|\\
&+&\|v(t)\phi\circ \af(a)- \phi\circ
\af(a)v(t)\|\\
&&\hspace{-0.5in}=\|\bt^2(v(t))\cdots \bt^{k-1}(v(t))\phi\circ \af(a)-
\phi\circ \af(a)\bt^2(v(t))\cdots\bt^{k-1}(v(t))\|\\
&+&\|\bt(v(t))\phi\circ \af(a)
-\phi\circ \af(a)\bt(v(t))\|\\
&+&\|v(t)\phi\circ \af(a)- \phi\circ
\af(a)v(t)\|.\\
\eneq
It follows that
\beq\label{LBM-8}
\|[u_k(t),\,\phi\circ \af(a)]\|\le
\sum_{j=0}^{k-1}\|\bt^j(u(t))\phi\circ \af(a)-\phi\circ
\af(a)\bt^j(u(t))\|.
\eneq
Thus, for $k\le N_1,$ by (\ref{LBM-6}) and by (\ref{LBM-2-2})
\beq\label{LBM-9}
\hspace{-0.8in}\|[u_k(t),\,\phi\circ \af(a)]\|&\le&
\sum_{j=0}^{k-1}\|[\bt^j(u(t)),\,\phi\circ \af(a)]\|\\
&&\hspace{-0.8in}\le \sum_{j=0}^{k-1}(\|\bt^j(v(t))\bt^j\circ\phi\circ \af^{-j+1}(a)-
\bt^j\circ \phi\circ \af^{-j+1}(a)\bt^j(v(t))\|+j\dt_1)\\
&&\hspace{-0.8in}\le \sum_{j=0}^{k-1}(\|v(t)\phi\circ \af^{-j+1}(a)-
\phi\circ \af^{-j+1}(a)v(t)\|+j\dt_1)\\
&&\hspace{-0.8in}\le N_1(\dt_1/2+(N_1-1)\dt_1)<N_1^2\dt_1<\eta_1
\eneq
for all $a\in {\cal G}_1$ and $t\in
[0,t_1].$

Applying the property (H) for $L>0,$  we obtain a continuous path
of unitaries $\{w(t): t\in [0,1]\}$ be such that
$$
w(0)=u_{N_1}(t_1),\,\,\,w(1)=1
$$
\beq\label{LBM-10}
\|w(t)-w(t')\|\le L|t-t'|\andeqn \|[w(t),\,\phi\circ
\af(a)]\|<{1\over{16}}
\eneq
for all $t, t'\in [0,1]$ and $a\in \cup_{i=1}^{N_1-1}\af^{-i}({\cal
F}_1).$

Let $0=s_0<s_1<\cdots <s_{N_1-1}=1$ such that
\beq\label{LBM-10+}
|s_i-s_{i-1}|={1\over{N_1-1}},\,\,\,i=1,2,...,N_1-1.
\eneq
Put $S=\{u_j(t_1), w(s_j): j=0,1,2,...,N_1-1\}.$
Since $\bt$ has the cyclic Rokhlin property (see \ref{DU}), we obtain mutually orthogonal
projections $e_1,e_2,...,e_{N_1}\in B$ such that
$\sum_{i=1}^{N_1}e_i=1_B$ and
\beq\label{LBM-11}
\|[e_i,\, b]\|<{L\over{128N_1(N_1-1)}},\,\,\,i=1,2,...,N_1
\eneq
for $b\in \cup_{i=0}^{N_1-1}\bt^i(S)\cup\cup_{j=0}^{N_1}\af^{-j}({\cal G}_1)\cup \phi(\cup_{j=0}^{N_1}\af^{-j+1}({\cal G}_1))$ and with $e_{N_1+1}=e_1,$
\beq\label{LBM-12}
\|\bt(e_i)-e_{i+1}\|<{L\over{128N_1(N_1-1)}},\,\,\,i=1,2,..., N_1.
\eneq
Define
\beq\label{LBM-13}
Z=\sum_{j=1}^{N_1}u_{j-1}(t_1)\bt^{j-1}(w(s_{j-1}))e_j.
\eneq

Denote $\ep'={L\over{128N_1(N_1-1)}}.$
We estimate that, using (\ref{LBM-12}) and (\ref{LBM-11})
and the fact that $e_ie_j=0$ (for $i\not=j$), as well as
(\ref{LBM-10}) and (\ref{LBM-10+}),
\beq\label{LBM-14}
\hspace{-0.4in}Z\bt(Z)^*&\approx_{3N_1\ep'}&\sum_{j=1}^{N_1} e_ju_{j-1}(t_1)\bt^{j-1}(w(s_{j-1}))e_j
\bt^{j-1}(w(s_{j-2})^*)\bt(u_{j-2}(t_1)^*)\\\label{LBM-15}
\hspace{-0.4in}&\approx_{2N_1\ep'+{L\over{N_1-1}}}&\sum_{j=1}^{N_1}e_ju_{j-1}(t_1)\bt(u_{j-2}(t_1)^*)e_j\\\label{LBM-16}
\hspace{-0.4in}&\approx_{N_1\ep'}&\sum_{j=1}^{N_1}v(t_1)e_j=v(t_1).
\eneq
Similarly, we then estimate that
\beq\label{LBM-17}
\|Z^*Z-1_B\|<{L\over{32(N_1-1)}}\andeqn \|ZZ^*-1_B\|<{L\over{32(N_1-1)}}.
\eneq
It follows that there is a unitary $V_1\in B$ such that
\beq\label{LBM-18}
\|V_1-Z\|<{L\over{16(N_1-1)}}.
\eneq
We also estimate that
\beq\label{LBM-15-}
Z\phi\circ \af(a) &=&\sum_{j=1}^{N_1}u_{j-1}(t_1)\bt^{j-1}(w(s_{j-1}))e_j\phi\circ \af(a)\\\label{LBM-15+1}
&\approx_{3N_1\ep'+{N_1\dt_1\over{2}}}&\sum_{j=1}^{N_1}e_ju_{j-1}(t_1)\bt^{j-1}(w(s_{j-1})\phi\circ \af^{-j+2}(a))e_j\\
\label{LBM-15+2}
&\approx_{1/16}&\sum_{j=1}^{N_1}e_ju_{j-1}(t_1)\bt^{j-1}\circ \phi\circ\af^{-j+2}(a)\bt^{j-1}(w(s_{j-1}))e_j\\\label{LBM-15+3}
&\approx_{N_1\dt_1\over{2}}&\sum_{j=1}^{N_1}e_ju_{j-1}(t_1)\phi\circ \af(a)\bt^{j-1}(w(s_{j-1}))e_j\\\label{LBM15+4}
&\approx_{\dt_1/2}& \sum_{j=1}^{N_1}e_j\phi\circ \af(a)u_{j-1}(t_1)\bt^{j-1}(w(s_{j-1}))e_j\\\label{LBM-15+4}
&\approx_{3N_1\ep'}&\phi\circ\af(a) Z
\eneq
for all $a\in {\cal F}_1.$

Note that, since $e_{N_1+1}=e_1,$
$$
u_0(t_1)w(0)e_1=u_{N_1}(t_1)e_1=u_{N_1}(t)\bt^{N_1}(w(t_{N_1-1}))e_{N_1+1}.
$$
Since  $\{e_i: 1\le i\le N_1\}$ is mutually orthogonal, as above,
we also have
\beq\label{LBM-20}
v(t_1)\bt(Z) &=&
v(t_1)\sum_{j=1}^{N_1}\bt(u_{j-1}(t_1))\bt^j(w(t_{j-1}))\bt(e_j)\\
&\approx_{4N_1\ep'}&\sum_{j=1}^{N_1}e_{j+1}v(t_1)\bt(u_{j-1}(t_1))\bt^j(w(t_{j-1}))e_{j+1}\\
&\approx_{{L\over{N_1-1}}}&\sum_{j=1}^{N_1}e_{j+1}v(t_1)\bt(u_{j-1}(t_1))\bt^j(w(t_j))e_{j+1}\\
&\approx_{2N_1\ep'}&\sum_{j=1}^{N_1}u_j(t_1)\bt^j(w(t_j))e_{j+1}\\
&=& Z.
\eneq
Therefore (by (\ref{LBM-18}))
\beq\label{LBM-21}
\|v(t_1)\bt(V_1)-V_1\|<({(4+128+2)L\over{128(N_1-1)}})+{2L\over{16(N_1-1)}}<{1\over{8}}.
\eneq

It follows from
that
\beq\label{LBM-22}
\|v(t_1)-V_1\bt(V_1)^*\|<1/8.
\eneq
By  (\ref{LBM-15-}), (\ref{LBM-15+1}), (\ref{LBM-15+2}),  (\ref{LBM-15+3}) , (\ref{LBM15+4}), (\ref{LBM-15+4}) and (\ref{LBM-18})
\beq\label{LBM-23}
\|[V_1, \phi\circ
\af(a)]\|&<&3N_1\ep'+{N_1\dt_1\over{2}}+{1\over{16}}+{N_1\dt_1\over{2}}+{\dt_1\over{2}}+
3N_1\ep'+{L\over{8(N_1-1)}}\\\label{LBM-23+}
&<&1/4
\eneq
for all $a\in {\cal F}_1.$

Define
\beq\nonumber
&&z_1=V_1^*v(t_1)\bt(V_1),\,\,\,\phi_1={\rm ad}\, V_1\circ \phi,\\\nonumber
&&\bt_1={\rm ad}\, z_1^*\circ \bt\andeqn
u_1(t)=V_1^*v(t_1)v(t+t_1)V_1.
\eneq
Then
\beq\nonumber
&&\hspace{-0.4in}{\rm ad}u_1(t)\circ  \bt_1\circ \phi_1(a) =
u_1(t)^*z_1(\bt(V_1)^*\bt\circ \phi(a)\bt(V_1))
z_1^*u_1(t)\\\nonumber
&=&V_1^*v(t+t_1)^*v(t_1)^*V_1 (V_1^*v(t_1)\bt\circ
\phi(a)v(t_1)^*V_1)V_1^*v(t_1)v(t+t_1)V_1\\\label{LBM-24}
&=& V_1^*v(t+t_1)^*\bt\circ \phi(a) v(t+t_1)V_1
 \eneq
for all $a\in A.$ It follows that
\beq\label{LBM-25}
\lim_{t\to\infty}{\rm ad} u_1(t)\circ \bt_1\circ \phi_1(a)&=&\lim_{t\to\infty} {\rm ad}\, (v(t+t_1)V_1)\circ \bt\circ \phi(a)\\\label{LBM-26}
&=&{\rm ad}\, V_1\circ \phi\circ \af(a)=\phi_1\circ \af(a)
\eneq
for all $a\in A.$

We also have $u_1(0)=1.$ By (\ref{LBM-21}),
\beq\label{LBM-27}
\|z_1-1\|=\|V_1^*v(t_1)\bt(V_1)-V_1^*V_1\|<{1\over{8}}.
\eneq

By  (\ref{LBM-3}),
\beq\nonumber
&&\|\phi_1\circ \af(a)-{\rm ad}\, u_1(t)\circ \bt_1\circ
\phi_1(a)\|\\\label{LBM-28}
&=&\|\phi\circ \af(a)-{\rm ad}\, v(t+t_1)\circ \bt\circ
\phi(a)\|<{\dt_2\over{2}}
\eneq
for all $a\in \cup_{i=0}^{N_2-1}({\cal G}_2)$ and all $t\in [0,\infty).$

With (\ref{LBM-25}), (\ref{LBM-26}) and (\ref{LBM-28}), we will
repeat this process. To repeat the process, we also note that, by
\ref{DU}, $\bt_1$ has the cyclic Rokhlin property in the sense of
\ref{DU}.

By repeating the process, we obtain sequences of unitaries
$\{V_n\}$  and $\{z_n\}$ of $B$ such that
\beq\label{LBM-29}
\|[V_n, \phi\circ \af(a)\|<{1\over{4^n}}
\eneq
for all $a\in {\cal F}_n,$
\beq\label{LBM-30}
\|z_n-1\|<{1\over{4^n}},
\eneq
$n=1,2,....$ Moreover, with
$$
\phi_n={\rm ad}\, (V_1V_2\cdots V_{n-1}V_n)\circ \phi\andeqn
\bt_n={\rm ad}\,(z_1^*z_2^*\cdots z_{n-1}^*z_n^*)\circ \bt,
$$
\beq\label{LBM-31}
\phi_n\circ \af\approx_{{\dt_n\over{2}}} \bt_n\circ
\phi_n\,\,\,\text{on}\,\,\, {\cal F}_n,
\eneq
$n=1,2,....$ By (\ref{LBM-29}), there exists a \hm\, $\phi': A\to
B$ such that
$$
\phi'(a)=\lim_{n\to\infty}\phi_n(a)\rforal a\in A.
$$
Since each $\phi_n$ is a monomorphism, so is $\phi'.$ By
(\ref{LBM-30}), we obtain a unitary $w\in B$ such that
\beq\label{LBM-32}
w=\lim_{n\to\infty}z_1^*z_2^*\cdots z_{n-1}^*z_n^*.
\eneq
It follows from (\ref{LBM-25})  and (\ref{LBM-26}) that
$$
\phi'\circ \af={\rm ad}\, w\circ \bt \circ \phi'.
$$

\end{proof}

\begin{cor}\label{MC}
Let $A$ be a unital AH-algebra, let $\af\in Aut(A)$ and let $B\in {\cal N}$ be
a unital separable simple \CA\, with tracial rank zero. Suppose that there is a unital
monomorphism $h: A\to B$ such that there exists a continuous path
of unitaries $\{U(t): t\in [0, \infty)\}$ of $B$ satisfying the
following
$$
\lim_{t\to\infty}{\rm ad}\, U(t)\circ h\circ \af(a)=h(a)\tforal
a\in A.
$$
Then $A\rtimes_{\af}\Z$ can be embedded into a unital simple
AF-algebra $C.$
\end{cor}

\begin{proof}
We first note that $A\otimes {\cal U}$ also satisfies the conditions in the statement above  with $\af$ replaced by $\af\otimes \sigma$ and $B$ replaced by $B\otimes {\cal U}, $  and $h$ by the  monomorphism 
$h'$ defined by $h'(a\otimes b)=h(a)\otimes b$ for all $a\in A$ and $b\in {\cal U}.$ 
So we may assume that $\af$ satisfy the cyclic Rokhlin property. 
Next we  note that, by \ref{path}, any unital monomorphism from
any  unital AH-algebra to a unital separable simple \CA\, with
tracial rank zero has  property (H) for $L=2\pi+1.$ So, in
particular, $h$ has the property (H). Let $\tau\in T(B)$ be a
tracial state and let $\D$ be the range of $\tau$ on $K_0(B).$
There is a unital simple AF-algebra $B_1$ with a unique tracial
state such that $(K_0(B_1), K_0(B_1)_+, [1_{B_1}])= (\D, \D_+,
1).$ By the classification theorem of \cite{Lnduke}, $B$ is a
unital simple AH-algebra with real rank zero and with no dimension
growth. By \cite{EG},
 there is a unital monomorphism $\phi: B\to
B_1.$ Denote by $C=B_1\otimes {\cal U}$ and
$\bt={\id}_{B_1}\otimes \sigma.$ It follows from \ref{LBM} 
(and using the cyclic Rokhlin property---see Lemma 4.1 of \cite{Lnemb2}) that
there is a unital monomorphism $\psi: A\rtimes_{\af}\Z\to
C\rtimes_{\bt}\Z.$ Since $(\bt)_{*0}={\id}_{K_0(C)}$ and $\bt$ has
the cyclic Rokhlin property, it follows from Theorem 3.4 of
\cite{Lngpots} that $C\rtimes_{\bt}\Z$ has tracial rank zero. It
follows from the classification theorem (see \cite{Lnduke} and
\cite{Lnann}) that $C\rtimes_{\bt}\Z$ is a unital simple
AH-algebra with slow dimension growth and real rank zero. It
follows from \cite{EG} that $C\rtimes_{\bt}\Z$ can be embedded
into a unital simple AF-algebra. It follows that
$A\rtimes_{\af}\Z$ can be embedded into a unital simple AF-algebra
(with a unique tracial state).

\end{proof}

\begin{lem}\label{kappa}
Let $A$ be a unital separable amenable \CA. Then for any finite subset
${\cal P}\subset K_1(A),$ there exists $\dt>0$ and a finite subset ${\cal G}\subset A$ satisfying the following:
For any pair of automorphisms $\af_1, \af_2\in Aut(A),$ any unital \hm\, $h : A\to B$ (where $B$ is a unital \CA)
and any unitary
$u\in U(B),$ if
\beq\label{kappa1}
\|[h\circ \af_i(a), \,u]\|<\dt\,\,\,\tforal a\in {\cal G}\,\,\,and\,\,\,
(\af_1)_{*1}=(-\af_2)_{*1},
\eneq
then
\beq\label{kappa2}
\rm{bott}_1(h\circ \af_1, u)|_{\cal
P}=-\rm{bott}_1(h\circ\af_2,u)|_{\cal P}.
\eneq
\end{lem}
(see 2.10 of \cite{Lnhomp} for the definition of
$\text{bott}_1(\phi, u)$)

\begin{proof}
Suppose that there is a unital \CA\, $B,$ there is a \hm\, ${\bar H}: A\to B$ and a unitary ${\bar U}\in B$ such that
$[{\bar H}(a),{\bar U}]=0$ for all $a\in A.$

Let ${\bar \Phi}: A\otimes C(S^1)\to B$ be the \hm\, defined by ${\bar \Phi}(a\otimes f)={\bar H}(a)f({\bar U})$ for all $a\in A$ and $f\in C(S^1).$
Let $\bt:K_1(A)\to K_0(A\otimes C(S^1))$ be the \hm\, given by the tensor product.
Then $\text{bott}_1({\bar H}\circ \af_i,{\bar U})={\bar \Phi}_{*0}\circ \bt \circ (\af_i)_{*1},$ $i=1,2.$
It follows that
\beq\label{kappa3}
\text{bott}_1({\bar H}\circ \af_1, {\bar U})=-\text{bott}_1({\bar H}\circ \af_1, {\bar U}).
\eneq

To see the lemma holds, we assume that it fails. We then obtain a
finite subset ${\cal P}\subset K_1(A),$ a sequence of unital \CA s
$B_n,$ a sequence of decreasing positive numbers $\{\dt_n\}$ with
$\lim_{n\to\infty}\dt_n=0,$ an increasing sequence of finite
subsets ${\cal G}_n\subset A$ such that $\cup_{n=1}^{\infty} {\cal
G}_n$ is dense in $A$ and a sequence of unital  \hm s $h_n: A\to
B_n,$ and a sequence of unitaries $u_n\in B_n$ such that
\beq\label{kappa4}
\|[h_n(a), u_n]\|<\dt_n\tforal a\in {\cal G}_n\andeqn
\eneq
\beq\label{kappa5}
\text{bott}_1(h_n\circ \af_1, u_n)|_{{\cal P}}\not=-\text{bott}_1(h_n\circ \af_2, u_n)|_{{\cal P}}
\eneq
(for all sufficiently large $n$).

Define a map $H: A\to l^{\infty}(\{B_n\})$ and a unitary $U=\{u_n\}.$ Denote by
$\pi: l^{\infty}(\{B_n\})\to q_{\infty}(\{B_n\})$ the quotient map.
Denote by ${\bar H}=\pi\circ H$ and ${\bar U}=\pi(U).$ Then
$$
[{\bar H}(a), {\bar U}]=0\tforal a\in A.
$$
It follows from what has been proved,
\beq\label{kappa6}
\text{bott}_1({\bar H}\circ\af_1, {\bar U})=-\text{bott}_1({\bar H}\circ \af_2, {\bar U}).
\eneq

Then, for any sufficiently large $n,$
\beq\label{kappa7}
\text{bott}_1(h_n\circ \af_1, u_n)|_{\cal P}=-\text{bott}_1(h_n\circ \af_2, u_n)|_{{\cal P}}.
\eneq
This contradicts with (\ref{kappa5}).

\end{proof}

\begin{lem}\label{AHembL1}
Let $A$ be a unital AH-algebra and let $B$ be a unital separable
simple \CA\, of tracial rank zero with $K_1(B)=\{0\}$ and
divisible $K_0(B).$ Let $\af, \kappa\in Aut(A)$ be  automorphisms
with $\kappa_{*1}(x)=-x$ for all $x\in K_1(A).$  Suppose that
there is a unital monomorphism $\phi: A\to B$ such that
\beq\label{AHembL-1}
[\phi\circ \af]=[\phi] \,\,\,\text{in}\,\,\, KL(A,B)
\eneq
and
\beq\label{AHembL-2}
\tau\circ \phi\circ \af=\tau\circ \phi
\eneq
for all $\tau\in T(B).$ Suppose also that there is a continuous
path of unitaries $\{V(t): t\in [0,\infty)\}$ of $B$ such that
\beq\label{AHembL-3}
\lim_{t\to\infty}{\rm ad}\, V(t)\circ \phi\circ \af\circ \kappa(a)
=\phi\circ \kappa \circ \af(a)
\eneq
for all $a\in A.$

 Then there is a unital monomorphism $\psi: A\to M_2(B\otimes
{\cal U})$ and there is a unitary $w\in U(M_2(B\otimes {\cal U}))$
such that
$$
\psi\circ \af={\rm ad}\, w\circ ({\rm id}_B\otimes \sigma)\circ
\psi.
$$
\end{lem}

\begin{proof}
To simply notation, we may assume that $B=B\otimes {\cal U}$ and
denote $\bt={\rm id}_B\otimes \sigma.$ It is well known that $\bt$
is approximately inner.

It follows from Theorem \ref{cdT1}, (\ref{AHembL-1}) and
(\ref{AHembL-2}) that there exists a sequence of unitaries
$\{u_n\}\subset B$ such that
\beq\label{AHeL-1}
\phi\circ \af(a)=\lim_{n\to\infty}{\rm ad}\, u_n\circ \bt\circ
\phi(a)
\eneq
for all $a\in A.$

Define $\psi_0, {\bar \psi_0}: A\to M_2(B)$ by
\beq\label{AHeL-2}
\psi_0(a)=\begin{pmatrix} \phi(a) & 0\\
                             0 & \phi\circ \kappa
                             \end{pmatrix}
                            \andeqn \eneq
  \beq\label{AHeL-2-1}
  {\bar \psi}_0=\begin{pmatrix} \phi\circ \af(a) & 0 \\
                                 0 & \phi\circ \af\circ \kappa(a)
                                 \end{pmatrix}
                                 \eneq
                             for all $a\in A.$
Define $\bt^{(2)}=\bt\otimes {\rm id}_{M_2}.$

It follows from (\ref{AHeL-1}), by passing to a subsequence if
necessary, with $U_n={\rm diag}(u_n, u_n),$ we may assume that
\beq\label{AHeL-3}
{\bar \psi}_0(a)=\lim_{n\to\infty}{\rm ad}\, U_n\circ
\bt^{(2)}\circ \psi_0(a)
\eneq
for all $a\in A.$

Let ${\cal F}_1\subset {\cal F}_2\cdots \subset {\cal F}_n\subset
\cdots $ be an increasing sequence of the unit ball of $B$ such
that the union is dense in the unit ball of $B.$

Let $1>\ep_n>0$ such that $\sum_{n=1}^{\infty}\ep_n<\infty.$  For
$\psi_0: A\to M_2(B),$ there exists $\dt_n'>0,$ a finite subset
${\cal G}_n'\subset A$  and a finite subset ${\cal P}_n\subset
\underline{K}(A)$ required by \ref{BHPT} for $\ep_n$ and ${\cal
F}_n.$ By choosing larger ${\cal P}_n,$ we may assume that if
$x\in {\cal P}_n,$ $-x\in {\cal P}_n.$
For each $n,$ choose $0<\dt_n\le \dt_n'$ and a finite subset ${\cal G}_n\supset {\cal G}_n'$ required
in \ref{kappa} for ${\cal P}_n\cap K_1(A).$

Without loss of generality, we may assume that ${\cal G}_n$ are in
the unit ball of $A,$ ${\cal F}_n\subset {\cal G}_n\subset {\cal
G}_{n+1},$ $n=1,2,....$ Put $S=\{1_{C(S^1)}, z\},$ where $z\in
C(S^1)$ is the identity function on the unit circle. Here we also
assume that for any $\dt_n$-${\cal G}_n\otimes S$-multiplicative
\morp s $L_1, L_2: A\otimes C(S^1)\to B$ (for any unital \CA\,
$B$), $[L_i]|_{{\cal P}_n}$ is well defined and
$$
[L_1]|_{{\cal P}_n}=[L_2]|_{{\cal P}_n},
$$
provided that
$$
L_1\approx_{\dt_n} L_2\,\,\,\text{on}\,\,\,{\cal G}_n.
$$

We may assume that
\beq\label{AHeL-4}
{\bar \psi}_0\approx_{\dt_n/4} {\rm ad}\, U_n\circ \bt^{(2)}\circ
\psi_0\,\,\,\text{on}\,\,\,{\cal G}_n.
\eneq
Put $W(t)={\rm diag}(1, V(t)).$  There are $t_n<t_{n+1}\to\infty$
such that
\beq\label{AHeL-5}
{\rm ad}\, W(t)\circ {\bar \psi}_0\approx_{\dt_n/4} \psi_0\circ
\af\,\,\, \text{on}\,\,\,{\cal G}_n
\eneq
for all $t\ge t_n,$ $n=1,2,....$

Put
\beq\label{AHeL-6}
Z_n=\begin{pmatrix} u_n^*u_{n+1} & 0\\
                      0 & u_n^*V(t_n)^*V(t_{n+1})u_{n+1}
                      \end{pmatrix}.
                      \eneq
and
\beq\label{AHeL-7}
Z_n(t)=\begin{pmatrix}u_n^*u_{n+1} & 0\\
                      0 &
                      u_n^*V(t_n)^*V(t_{n+1}+(t_n-t_{n+1})t)u_{n+1}
                      \end{pmatrix}.
                      \eneq
for $t\in [0,1].$ So $\{Z_n(t): t\in [0,1]\}$ is a path of
continuous unitaries, $Z_n(0)=Z_n$ and
\beq\label{AHeL-8}
Z_n(1)=\begin{pmatrix}u_n^*u_{n+1} & 0\\
                      0 & u_n^*u_{n+1}
                      \end{pmatrix}=U_n^*U_{n+1}.
                      \eneq
Note that
\beq\label{AHeL-9}
\|[\psi_0\circ \af(a),\, Z_n(t)]\|<\dt_n\rforal a\in {\cal
G}_n\andeqn t\in [0,1].
\eneq
It follows that
\beq\label{AHeL-10}
\text{Bott}(\psi_0, Z_n)|_{{\cal P}_n}=\text{Bott}(\psi_0,
U_n^*U_{n+1})|_{{\cal P}_n}.
\eneq
By the choice of $\dt_n$ and ${\cal G}_n,$ by applying \ref{kappa},
\beq\label{AHeL-11}
\text{bott}_1(\psi_0, Z_n)|_{{\cal P}_n}&=&\text{bott}_1(\psi_0,
U_n^*U_{n+1})|_{{\cal P}_n}\\
&=&\text{bott}_1(\phi\circ \af, u_n^*u_{n+1})+
\text{bott}_1(\phi\circ \kappa\circ \af, u_n^*u_{n+1})\\
&=& 0.
\eneq
Since $K_1(B)=\{0\}$ and $K_0(B)$ is divisible, the above implies
that
\beq\label{AHeL-12}
\text{Bott}(\psi_0\circ \af, Z_n)|_{{\cal P}_n}=0,\,\,\,n=1,2,....
\eneq
It follows \ref{BHPT} that there is a rectifiable continuous path
of unitaries $\{Z_n'(t): t\in [0,1]\}$ such that
\beq\label{AHeL-13}
Z_n'(1)=Z_n,\,\,\,Z_n'(0)=1\andeqn \|[\psi_0(a), Z_n'(t)]\|<\ep_n
\eneq
for all $a\in {\cal F}_n$ and $t\in [0,1].$ Define
\beq\label{AHeL-14}
v(t)=U_nZ_n'(t-n)\rforal t\in [n, n+1],\,\,\, n=1,2,....
\eneq
Note that $v(n)=U_nZ_n'(0)=U_n,$ $v(n+1)=U_nZ_n'(1)=U_{n+1}$ and
$V(t)$ is continuous on $[1, \infty).$ We also compute that
\beq\label{AHeL-15}
\lim_{t\to\infty} {\rm ad}\, v(t)\circ \bt\circ
\psi_0(a)=\psi_0\circ \af(a)
\eneq
for all $a\in A.$

It follows from \ref{path} that $\psi_0\circ \af$ has the property
(H) for the constant $2\pi+1.$ Thus Lemma \ref{LBM} applies.

\end{proof}

\begin{thm}\label{AHembT1}
Let $A$ be a unital AH-algebra (not necessarily simple) and let
$\af\in Aut(A)$ be a automorphism. Suppose that there is another
automorphism $\kappa\in Aut(A)$ such that $\kappa(x)=-x$ for all
$x\in K_1(A)$ such that there is a continuous path of unitaries
$\{U(t): t\in [0, \infty)\}$ such that
$$
\lim_{t\to\infty}{\rm ad}\, U(t)\circ \af\circ
\kappa(a)=\kappa\circ \af(a)
$$
for all $a\in A.$ Then $A\rtimes_{\af}\Z$ can be embedded into a
unital simple AF-algebra if and only if $A$ has a faithful
$\af$-invariant tracial state.
\end{thm}

\begin{proof}
Suppose that $\phi: A\rtimes_{\af}\Z\to B$ is a unital
monomorphism, where $B$ is a unital simple AF-algebra. Let
$\tau\in T(B).$ Then $\tau\circ \phi$ gives a faithful
$\af$-invariant tracial state on $A.$

 So now
suppose that there is an $\af$-invariant faithful tracial state
$\tau$ on $A.$

It follows from \ref{CembAH} that there is a unital simple
AF-algebra $B$ with a unique tracial state $\tau$ and $(K_0(B),
K_0(B)_+, [1_B])=(\D, \D_+, 1),$ where $\D$ is a dense additive
subgroup of $\R$ with $\D=\rho_B(K_0(B)),$ such that is a
unital embedding $\phi: A\to B$ with
\beq\label{AHeT-1}
\tau\circ \phi\circ \af=\tau\circ \phi.
\eneq
By considering $B\otimes {\cal U},$ we may assume that $K_0(B)$ is
divisible. It follows from that fact that $K_0(B)=\D$ and
(\ref{AHeT-1}),
\beq\label{AHeT-12}
(\phi\circ \af)_{*0}=\phi_{*0}.
\eneq
We also have $(\phi\circ \af)_{*1}=\phi_{*1}$ since
$K_1(B)=\{0\}.$
Since $K_0(B)$ is divisible, it follows that
$$
[\phi\circ \af]=[\phi]\,\,\,\text{in}\,\,\,KL(A, B).
$$

It follows from \ref{AHembL1}, there is a unital embedding $\psi:
A\to M_2(B\otimes {\cal U})$ and a unitary $w\in U(M_2(B\otimes
{\cal U}))$ such that
\beq\label{AHeT-3}
\phi\circ \af={\rm ad}\, w\circ ({\rm id}_B\otimes \sigma)\circ
\phi.
\eneq
Denote by $\bt={\rm id}_B\otimes \sigma$ and $C_1=M_2(B\otimes
{\cal U}).$ Then we have a unital embedding $\psi_1:
A\rtimes_{\af}\Z\to C_1\rtimes_{\bt}\Z.$ As in the proof of
\ref{MC}, this implies that $A\rtimes_{\af}\Z$ can be embedded
into a unital simple AF-algebra.


\end{proof}

It is worth to note that when $K_1(A)=\{0\},$ one can always take $\kappa={\rm id}_A.$

%

\begin{cor}\label{AHemC1}
Let $A$ be a unital AH-algebra and let $\af\in Aut(A).$ Suppose
there is $\kappa\in Aut(A)$ such that $\af(x)=-x$ for all $x\in
K_1(A)$ and $\af\circ \kappa=\kappa\circ \af.$  Then
$A\rtimes_{\af}\Z$ can be embedded into a unital simple AF-algebra
if and only if there is a faithful $\af$-invariant tracial state
\end{cor}
\begin{proof}
This immediately follows from \ref{AHembT1}.
\end{proof}

\begin{cor}\label{AHemC2}
Let $A$ be a unital AH-algebra and let $\af\in Aut(A).$ Suppose
that $\af(x)=-x$ for all $x\in K_1(A).$ Then $A\rtimes_{\af}\Z$
can be embedded into a unital simple AF-algebra if and only if
there is a faithful $\af$-invariant tracial state.

\end{cor}

\begin{proof}
Let $\kappa=\af.$ Then $\kappa\circ \af=\af\circ \kappa.$ So
\ref{AHemC1} applies.

\end{proof}

\begin{cor}\label{AHembC3}
Let $A$ be a unital AH-algebra (not necessarily simple) with
torsion $K_1(A).$ Let $\af\in Aut(A).$ Then $A\rtimes_{\af}\Z$ can
be embedded into a unital simple AF-algebra if and only if $A$
admits a faithful $\af$-invariant tracial state.
\end{cor}

\begin{proof}
By embedding $A\rtimes_{\af}\Z$ into $A\otimes {\cal U}\rtimes
_{\af\otimes \sigma}\Z,$  we may assume that $K_1(A)=\{0\}.$ We
then choose $\kappa={\rm id}_A.$ So, either \ref{AHemC1} or
\ref{AHembT1} applies.
\end{proof}

\begin{cor}\label{AHembC4}
Let $A$ be a unital AF-algebra and let $\af\in Aut(A).$ Then
$A\rtimes_{\af}\Z$ can be embedded into a unital simple AF-algebra
if and only if $A$ admits a faithful $\af$-invariant tracial
state.
\end{cor}

\begin{proof}
Every AF-algebra is an AH-algebra with trivial $K_1(A).$
\end{proof}

\begin{rem}
{\rm

The scheme of the proof of this section is  the same as that in
\cite{M1}.
 However, we use \ref{BHPT} and \ref{path} as well as \ref{cdT1} and
 \ref{CembAH}. Moreover, by
introducing $\kappa,$ we are able to avoid the difficult
generalization of results in \cite{KK}. One can verify, for
example, for unital simple A$\T$-algebra, required $\kappa$ can be
produced. However, in the next section, we will show that if $A$
is a unital A$\T$-algebra and $\af\in Aut(A)$ is an automorphism
then $A\rtimes_{\af}\Z$ can always be embedded into a unital
simple AF-algebra provided it has a faithful $\af$-invariant
tracial state.

 It probably worth to point out that, even simple
AH-algebras with trivial $K_1$-groups and torsion free
$K_0$-groups are   not necessary AF (see \cite{Li} for a
classification for some of those simple AH-algebras). In
particular, they may not have real rank zero.  There are also
simple AH-algebras whose stable rank is higher than one. Some more
exotic AH-algebras were known. In general, when \CA s are not
assumed
 to be simple,
there is no realistic hope of classification by $K$-theory at this moment.

}
\end{rem}

\section{Crossed products of A$\T$-algebras by $\Z$}

Some versions of the  following three facts are known. We collect here as
lemmas for the convenience.

\begin{lem}\label{group}
Let $A$ be a unital \CA\, and let $\tau\in T(A).$ Suppose that $u,
v\in U(A)$ are unitaries such that
\beq\label{gr1}
\|u-1\|<2/3,\,\,\,\|v-1\|<2/3\,\,\,\text{and}\,\,\,\|uv-1\|<2/3.
\eneq
Then
\beq\label{gr2}
\tau(\log(uv))=\tau(\log(u))+\tau(\log(v)).
\eneq
\end{lem}

\begin{proof}
Put
$$
\lambda=\max\{\|u-1\|,\|v-1\|, \|uv-1\|\}.
$$

Let
$$
F=\{e^{\pi i t}:-\pi+2\theta\le t\le \pi-2\theta\},
$$
where $\theta=\arccos (3\lambda/2).$

Then $F$  contains spectra of $u,$ $v$ and $uv.$ Take a branch of
logarithm $\log: F\to (-\pi, \pi)$ which is continuous (smooth) on $F.$
Denote
$$
b_1={1\over{2\pi i}}\log(u),\,\,\, b_2={1\over{2\pi
i}}\log(v)\andeqn b_3={1\over{2\pi i}}\log (uv).
$$
Set
$$
w(t)=\begin{cases} e^{-6\pi i t b_1} & 0\le t\le 1/3\\
                   e^{-2\pi i b_1} e^{-2\pi i (3t-1) b_2} & 1/3<t\le 2/3\\
                   e^{-2\pi i b_1} e^{-2\pi i  b_2}e^{2\pi i(3t-2)b_3} &
                   2/3<t\le 1.
                   \end{cases}
$$
Note that, for each $t\in [0,1],$ ${\rm sp}(w(t))\subset F.$ So
$\log w(t)$ is smooth and $ \log(w(0))=\log(w(1)).$
One then computes
\beq\label{gr-1}
{1\over{2\pi i}}\int_0^1 \tau( {d(w(t))\over{dt}} w(t)^*)dt=0.
\eneq

However, one also has that
\beq\label{trhp3}
&&\hspace{-0.3in}{1\over{2\pi i}}\int_0^1 \tau( {d(w(t))\over{dt}} w(t)^*)dt
={1\over{2 \pi i}}\{\int_0^{1/3} -6\pi i \tau(b_1) dt +\\
&+&\int_{1/3}^{2/3} 6\pi i \tau(b_2) dt
+\int_{2/3}^1 6\pi i \tau(b_3)  dt\}\\
&=&-\tau(b_1)-\tau(b_2)+\tau(b_3).
\eneq

\end{proof}

\begin{lem}\label{tracehp0}
 Let $A$ be a unital \CA, let $B\subset
A$ be a unital \SCA\, (with $1_B=1_A$),  let $\af: B\to B$ be a
unital \hm\, and let $\tau\in T(A)$ be a tracial state such that
$\tau(\af(a))=a$ for all $a\in B.$ Suppose that $b\in B_{s.a}$
with $\|b\|<{\arcsin(1/6)\over{\pi}},$  $w\in U(B)$  and $u\in
U(A)$ such that
\beq\label{trhp0}
\|w{\rm ad} u(\af(w^*))-1\|<1/3\andeqn  \|e^{2\pi i b}w{\rm ad}\,
u(\af(w^*)e^{-2\pi i\af(b)})-1\|<2/3.
\eneq

Then
\beq\label{trhp1}
\tau(\log(e^{2\pi i b}w{\rm ad}\, u(\af((e^{-2\pi ib}w)^*)))=
\tau(\log(w{\rm ad}\, u(\af(w^*)))).
\eneq
\end{lem}

\begin{proof}
By applying \ref{group}, we have
\beq\label{trhp0-1}
&&\hspace{-0.4in}\tau(\log(e^{2\pi i b}w{\rm ad}\, u(\af(e^{2\pi
i b} w)^*)))\\
&=& \tau(\log(e^{2\pi i b}))+
\tau(\log(w{\rm ad}\, u(\af(w^*)e^{-2\pi i\af(b)})))\\
&=& 2\pi i\tau(b)+\tau(\log(w{\rm ad}\, u(\af(w^*))
+\tau(\log({\rm ad}\, u(e^{-2\pi i\af(b)})))\\
&=& 2\pi i\tau(b)+\tau(\log(w{\rm ad}\, u(\af(w^*))
-2\pi i \tau(\af(b))\\
&=& \tau(\log(w{\rm ad}\, u(\af(w^*)).
\eneq

\end{proof}

\begin{lem}\label{tracehp2}
For each integer $k>0,$ there exists $\dt>0$ satisfying the
following:

Let $A$ be a unital \CA, let $B\subset A$ be a unital \SCA\, (with
$1_B=1_A$),  let $\af: B\to A$ be a unital \hm\, and let $\tau\in
T(A)$ be a  tracial state such that $\tau(\af(a))=a$ for all $a\in B,$
and let $z\in U(B).$

Suppose that there exists a unitary $u\in U(A)$ such that
\beq\label{trhp2-1}
\|{\rm ad}\, u\circ \af(z)-z\|<\dt.
\eneq
Then
\beq\label{trhp2-2}
\tau(\log(z^{m}{\rm ad}\, u(\af(z^{m*})))=m \tau(\log( z{\rm ad}\,
u(\af(z^*))))
\eneq
for all $|m|\le k.$

\end{lem}

\begin{proof}
 With sufficiently small $\dt,$
(\ref{trhp2-1}) implies that
\beq\label{trhp2-3}
\|z^m{\rm ad}\,u\circ \af(z^{m*})-1\|<1/2^{m}\andeqn \|{\rm
ad}\,u\circ \af(z^{m*})z^m-1\|<1/2^{m}
\eneq
for all $ |m|\le k.$
 Then, by \ref{group},
\beq\nonumber
\tau(\log(z^m{\rm ad}\,\af(z^{m*})))&=& \tau(\log(z^{m-1}{\rm
ad}\, u(\af(z^{m-1})^*){\rm ad}\, u(\af(z^*))z))\\\nonumber
&=&\tau(\log(z^{m-1}{\rm
ad}\, u(\af(z^{m-1})^*)))+\tau(\log({\rm ad}\, u(\af(z^*))z))\\
&=& \tau(\log(z^{m-1}{\rm
ad}\, u(\af(z^{m-1})^*)))+\tau(\log(z{\rm ad}\, u(\af(z^*)))).
\eneq
The lemma then follows by induction.

\end{proof}

Let $B=\oplus_{i=1}^N M_{r(i)}(C(X_i)),$ where $X_i$ is a compact
subset of the unit circle. Let $e_i\in M_{r(i)}(C(X_i))$ be a
constant rank one projection. Denote by $z_i$ the identity
function on $X_i.$ We identify $C(X_i)$ with
$e_iM_{r(i)}(C(X_i))e_i$ and $z_i$ with the unitary in
$e_iM_{r(i)}(C(X_i))e_i,$ $i=1,2,...,N.$ Denote by
$z_i'=(1-e_i)+z_i.$ Denote by $E_i$ the identity of $M_{r(i)},$
the $i$-th summand of $B.$

Suppose that $A$ is a unital \CA\, of stable rank one and
$B\subset A$ with $1_B=1_A.$ Let $\af: B\to A$ be a unital
monomorphism. Suppose that $u\in U(A)$  with $uE_i=E_iu,$
$i=1,2,...,N.$ Suppose that $\tau\in T(A)$ such that
$$
\tau\circ \af|_B=\tau|_B.
$$

 Let $Z=\{v_1,v_2,...,v_m\}$ be unitaries in $B.$ Then we have
the following:

\begin{prop}\label{PAT}
There exists $\dt>0$ and a finite subset ${\cal G}\subset B$
satisfying the following:

If
\beq\label{pat-1}
&&\|{\rm ad}\, u\circ \af(a)-a\|<\dt \rforal a\in {\cal G}
\andeqn \\\label{pat-1+}
&&\tau(\log(z_i{\rm ad}\, u(\af(z_i^*))))=0,\,\,\,i=1,2,...,N,
\eneq
then
\beq\label{pat-2}
\tau(\log(v_j{\rm ad}\,u(\af(v_j^*))))=0
\eneq
for all $v_j\in Z.$
\end{prop}

\begin{proof}
It is easy to see that we may assume that $E_i=1_A=1_B$ and $B$
has only one summand, $B=M_{r}(C(X)),$ where $X$ is a compact
subset of the unit circle. Let $e$ be a constant rank one
projection in $B.$ We now use $z$ for the identity function on $X$
and the corresponding unitary in $eBe$ and denote $z'=(1-e)+z.$

For each $v_j,$ we assume that $[v_j]=[z^{m_j}]$ in $K_1(B)$ for
some integer $m_j.$ There is a continuous path of unitaries
$\{v_j(t): t\in [0,1]\}$ such that
$$
v_j(0)=v_j\andeqn v_j(1)=z^{m_j},\,\,\,j=1,2,...,N.
$$
There are $0=t_0<t_1<\cdots < t_K=1$ such that
\beq\label{pat-3}
\|v_j(t_s)-v_j(t_{s-1})\|<1/3,\,\,\,\,s=1,2,...,K
\eneq
for $1\le j\le N.$

Let $\dt_1>0$ be in the lemma \ref{tracehp2} for $k=\max\{|m_j|:
j=1,2,...,N\}.$

Put
$$
\dt=\min\{\dt_1, {1\over{6^{K+1}}}\}.
$$
We now assume that ${\cal G}$ contains at least $z, z',$
$v_j(t_{s-1})^*v_j(t_s)$ and $v_j(t_s),$ $s=1,2,...,N$ and
$j=1,2,...,m.$ Now suppose that (\ref{pat-1}) and (\ref{pat-1+})
hold. In particular, we have
\beq\label{pat-4}
\|{\rm ad}\, u\circ \af(v_j(t_s))-v_j(t_s)\|<1/6^{K+1},
\eneq
$s=0,2,...,K$ and $j=1,2,...,N.$ Note that there is $h(j,s)\in
B_{s.a}$ with $\|h(j,s)\|<{\arcsin(1/6)\over{\pi}}$ such that
$v_j(t_{s-1})^*v_j(t_s)=e^{2\pi i h(j,s)}.$

By (\ref{pat-1+}), we have
\beq
\tau(\log(z'{\rm ad}\, u(\af(z')^*)))=0.
\eneq

It follows from \ref{tracehp2} that
\beq\label{pat-5}
\tau(\log((z')^{m_j}{\rm ad}\,
u(\af((z')^{m_j})^*)))=0,\,\,\,j=1,2,...,N.
\eneq

By applying \ref{tracehp0}, we obtain
\beq\label{pat-6}
\tau(\log(v_j(t_{N-1}){\rm ad}\, u\af(v_j(t_{N-1}))))=0.
\eneq
 Then the desired conclusion follows
the repeated application of \ref{tracehp0}.

\end{proof}

The following Exel's formula actually holds in much great
generality.

\begin{lem}\label{TRbot}
There is a positive number $\dt$ satisfying the following:

Let $A$ be a unital AF-algebra.
 Suppose that $u, v\in U(A)$ are
unitaries and suppose that
\beq\label{Trbot1}
\|uv-vu\|<\dt.
\eneq
Then $\rm{bott}_1(u,v)$ is well defined and
\beq\label{Trbot2}
\rho_A(\rm{bott}_1(u,v))(\tau)= {1\over{2\pi
i}}\tau(log(vuv^*u^*))\tforal \tau\in T(A).
\eneq

\end{lem}

\begin{proof}
There is $\dt_1>0$ such that $\text{bott}_1(u,v)$ is well defined
for any pair of unitaries $u$ and $v$ such that $\|uv-vu\|<\dt_1.$
Moreover, there are two projections $p, q\in M_2(A)$ such that
\beq\label{Trbot3}
\text{bott}_1(u,v)=[p]-[q].
\eneq

 Now we
suppose that $\|uv-vu\|<\dt_1/2.$ There is $\dt_2>0$ such that
\beq\label{Trbot4}
\text{bott}_1(u', v)=[p]-[q]
\eneq
for any other unitary $u'\in U(A)$ satisfying  $\|u'-u\|<\dt_2.$
Similarly, we obtain $\dt_3>0$ such that
\beq\label{Trbot4+1}
\text{bott}_1(u', v')=[p]-[q]
\eneq
for any pair of  unitaries $u', v'\in U(A)$ satisfying
$\|u'-u\|<\dt_3$ and $\|v'-v\|<\dt_3.$

 Now choose $\dt=\dt_1/2.$

 There are sequences  $\{u_n\}$  and $\{v_n\}$ of unitaries in
$\cup_{n=1}^{\infty}A_n$ such that
\beq\label{trbot5}
\lim_{n\to\infty}\|u-u_n\|=0\andeqn \lim_{n\to\infty}\|v-v_n\|=0.
\eneq
Without loss of generality, we may assume that
\beq\label{trbot6}
u_n, v_n\in A_n,\,\,\, \|u-u_n\|<\dt_3\andeqn
\|v-v_n\|<\dt_3,\,\,\,n\ge N_1
\eneq
for some sufficiently large $N_1.$
Then
\beq\label{trbot6+}
\text{bott}_1(u_n, v_n)=\text{bott}_1(u,v)=[p]-[q]\rforal n\ge N_1.
\eneq
 By choosing larger $n,$ without loss of
generality, we may assume that $p, q\in M_2(A_n).$ We may write
$A_n=\oplus_{k=1}^N M_{r(k)},$ where $r(k)\ge 1$ is an integer. We
may also write
\beq\label{trbot7}
p=\oplus_{k=1}^Np_k\andeqn q=\oplus_{k=1}^N q_k,
\eneq
where $p_k, q_k\in M_{r(k)}$ are projections, $k=1,2,...,N.$
Moreover,
\beq\label{trbot8}
u_n=\oplus_{k=1}^N u(k)\andeqn v=\oplus_{k=1}^N v(k),
\eneq
where $u(k),v(k)\in M_{r(k)}$ are unitaries,$k=1,2,...,N.$

 If follows
from Lemma 3.1 of \cite{Ex} that
\beq\label{trbot9}
Tr_k(p_k)-Tr_k(q_k)={1\over{2\pi
i}}Tr_k(\log(v(k)u(k)v(k)^*u(k)^*)),
\eneq
where $Tr_k$ denotes the standard trace on $M_{r(k)}.$ For any
tracial state $\tau\in T(A),$ there are non-negative numbers
$\lambda_1, \lambda_2,...,\lambda_{R(n)}$ with $ \sum_{k=1}^{R(n)}
{\lambda_k\over{r(k)}}=1 $ such that
\beq\label{trbot9+1}
\tau|_{A_n}=\sum_{k=1}^{R(n)} {\lambda_k Tr_k\over{r(k)}}.
\eneq
It follows that
\beq\label{trbot10}
\rho_A([p]-[q])(\tau)={1\over{2\pi i}}\tau(\log(v_nu_nv_n^*u_n^*))
\eneq
for each $n$ and all $\tau\in T(A).$  Note the left side remains
the same value.
 Since we choose a continuous branch of $\log,$  we conclude that
 \beq\label{trbot11}
 \rho_A(\text{bott}_1(u,v))={1\over{2\pi i}}\tau(\log(vuv^*u^*)).
 \eneq

\end{proof}

\begin{lem}\label{Asmu}
Let $B$ be a unital A$\T$-algebra and let $A$ be a unital simple
AF-algebra with a unique tracial state $\tau.$ Suppose that
$\af\in Aut(A)$ and suppose that there is a unital monomorphism
$\phi: B\to A$ such that
\beq\label{Asmu1}
\tau\circ \phi\circ \af=\tau\circ \phi.
\eneq
Then there is a unital simple AF-algebra $C$ with a unique tracial
state and a unital monomorphism $\psi: A\to C$ such that there
exists a continuous path of unitaries $\{U(t): t\in [0,\infty)\}$
of $C$ such that
\beq\label{Asmu2}
\lim_{t\to\infty}{\rm ad}\, U(t)\circ \psi\circ \phi\circ \af(a)=
\psi\circ \phi(a)
\eneq
for all $a\in B.$
\end{lem}

\begin{proof}
By replacing $A$ by $A\otimes {\cal U},$ we may assume that
$K_0(A)$ is divisible. Suppose that $\rho_A(K_0(A))=\D,$ where
$\D$ is a countable divisible dense subgroup of $\R.$ By embedding
$A$ into another simple AF-algebra, if necessary, we may assume
that $(K_0(A), K_0(A)_+, [1_A]) =(\D, \D_+, 1).$ In particular,
${\rm ker}\rho_A=\{0\}.$

Let $B={\overline{\cup_{n=1}^{\infty}B_n}},$ where
$$
B_n=\oplus_{i=1}^{R(n)} M_{r(i)}(C(X_{n,i}))
$$
and $X_{n,i}$ is a compact subset of the unit circle.

 Let $e^{(n,i)}$ be a constant rank one projection in
$M_{r(i)}(C(X_{n,i})).$ We will identify $C(X_{n,i})$ with
$e^{(n,i)}M_{r(i)}(C(X_{n,i}))e^{(n,i)}$ whenever it is
convenient. Denote by $z_{n,i}$ the identity function on $X_{n,i}$
and will be identified with   a unitary in
$e^{(n,i)}M_{r(i)}(C(X_{n,i}))e^{(n,i)},$ $i=1,2,...,R(n),$
$n=1,2,....$ Denote $z_{n,i}'=(1-e_{n,i})+z_{n,i}.$ To simplify
the notation, we will identify $B$ with $\phi(B).$ Also denote by
$D_n=\oplus_{i=1}^{R(n)} M_{r(i)}$ which is also identified with
\SCA\, of $B_n$ of constant matrices.

Since $K_1(A)=\{0\},$ both $A$ and $B$ have torsion $K$-theory,
${\rm ker}\rho_A=\{0\},$ by (\ref{Asmu1}), in addition to
(\ref{Asmu1}), we also have
\beq\label{asmu-1-1}
[\phi\circ \af]=[\phi]\,\,\, \text{in}\,\,\, KL(B, A).
\eneq
So combining (\ref{asmu-1-1}) with (\ref{Asmu1}), by \ref{cdT1},
there exists a sequence of unitaries $\{w_n\}\subset A$ such that
\beq\label{asmu-1}
\lim_{n\to\infty}{\rm ad}\, w_n \circ \af(b)=b
\eneq
for all $b\in B$ (note we have omitted the notation $\phi$).

Write $A=\overline{\cup_{n=1}^{\infty} A_n},$ where each $A_n$ is
a finite dimensional \SCA\, with $1_A\in A_n$ for each $n.$ We may
assume that $\{w_n\}\subset \cup_{n=1}^{\infty} A_n.$ Moreover, we
may assume that
\beq\label{asum-2}
{\rm ad}\, w_n\circ \af|_{D_n}={\id}_{D_n},\,\,\,n=1,2,....
\eneq

By passing to a subsequence if necessarily, we may also assume that
\beq\label{asum-3}
\|{\rm ad}\, w_{n+k}\circ \af(z_{n,i})-z_{n,i}\|
\eneq
is sufficiently small for all $k=1,2,...$ so that
$$
b_{n+k,i}={1\over{2\pi i}}\log
(z_{n,i}w_{n+k}^*\af(z_{n,i})w_{n+k})
$$
can be defined.

Let $\D_1$ be a countable subgroup of $\R$ which contains $\D$ and
$\{\tau(b_{n+k,i})\}.$ Let $C$ be a unital simple AF-algebra with
unique tracial state which we again denote it by $\tau$ such that
$(K_0(C), K_0(C)_+, [1_C])=(\D_1, (\D_1)_+,1).$ For the
convenience, we may also assume that $A\subset C$ (and $\tau(a)$
remains the same value) and we use $j: B\to C$ for the embedding.

Let $\{{\cal F}_n\}\subset B$ be an increasing sequence of finite
subsets whose union is dense in $B.$ Let $\ep_n>0$ (in place of
$\dt$), ${\cal G}_n$ (in place of ${\cal G}$) be a finite subset
of $B$ (in place of $A$)  and ${\cal P}_n\subset \underline{K}(B)$
be required in \ref{BHPT} for $1/n$ (in place of $\ep$) and ${\cal
F}_n$ (in place of ${\cal F}$) (and for $j$ in place of $h$
there). Since $B$ is an A$\T$-algebra, $K_i(B)$ is torsion free
($i=0,1$). Therefore, by the fact that $K_1(A)=\{0\}$ and $K_0(A)$ is divisible, it suffices to choose ${\cal P}_n\subset K_1(B).$ We
may assume that $\ep_n\downarrow 0.$ Without loss of generality,
we may also assume that ${\cal F}_n\subset B_n$ and ${\cal G}_n$ is in
the unit ball of $B_n.$ We may further assume that $z_{n,i}\in {\cal
G}_n.$

For each $n\ge 2,$ let $Z_n=\{z_{n,i}',z_{n-1,j}', i,j\}.$ To
simplify the notation, without loss of generality, we may assume
that $\ep_n$ is smaller than $\dt$ in \ref{PAT} and ${\cal G}_n$
is larger that ${\cal G}$ in \ref{PAT} corresponding $B_n$ and
$Z_n.$ Put $\gamma_n=\min\{\tau(e_{n,i}): 1\le i\le R(n)\}.$

There is $x_{n+k, i}\in K_0(C)$ such that
\beq\label{asum-4}
\tau(b_{n+k,i})=\rho_C(x_{n+k,i}), \,\,\,i=1,2,...,R(n), n,
k=1,2,....
\eneq

Let ${\cal G}_{n,i}\subset C(X_{n,i})$ be a finite subset such
that  ${\cal G}_n\subset \{\oplus_{i=1}^{R(n)}(a_{k,j}):
a_{k,j}\in {\cal G}_{n,i}\}.$ Note each $a_{k,j}$ may be viewed as
a function in $C(X_{n,i}).$

 Let $\dt_n'>0$ be required by Lemma 2 of \cite{M1} for
${\ep_n\over{2R(n)^2}}.$ In fact that we choose $\dt_n'$ even smaller so
that in the conclusion of Lemma 2 of \cite{M1} we may require that
\beq\label{asum-4+1+0}
\|[a_{k,j}(z),\,w]\|<\ep_n/2R(n)^2
\eneq
for all $a_{k,j}\in {\cal G}_{n,i}$ (in the case that
$X_{n,i}=S^1$).

Put $\dt_n=\min\{\ep_n/2, \gamma_n\cdot\dt_n',  2\sin (2\pi
\gamma_n\cdot \dt_n')\}.$

We also assume that the value in (\ref{asum-3}) is smaller than
$\dt_n.$

 We will construct a sequence of unitaries $\{u_n\}$ as
follows. For each $X_{n,i}\not=S^1,$  choose $\dt(n,i)>0$ so that
\beq\label{asum-4-1}
\text{bott}_1(z_{n,i}, z)=0
\eneq
for any unitary $z\in e_{n,i}Ce_{n,i},$  provided that
$$
\|[z_{n,i},\, z]\|<\dt(n,i).
$$
We then choose $w_{k(n)}$ so that
\beq\label{asum-4+1}
\|{\rm ad}\, w_{k(n)}\circ \af(z_{n,i})-z_{n,i}\|<\min\{\dt_n,
\dt(n,i)\}\\\label{asum-4+2}
 \andeqn\|{\rm ad}\, w_{k(n)}\circ
\af(a)-a\|<\dt_n\tforal a\in {\cal G}_n.
\eneq
Moreover,
\beq\label{asum-4+3}
\tau(b_{k(n),i})<\gamma_n\cdot \dt_n'.
\eneq

It follows from Lemma 2 of \cite{M1} that, if ${\rm
sp}(z_{n,i})=S^1,$  there is a unitary $v_{n,i}\in
e^{(n,i)}Ce^{(n,i)}$ such that
\beq\label{asum-5}
\|[z_{n,i}, v_{n,i}]\|<{\ep_n\over{R(n)^2}}\andeqn
\text{bott}_1(z_{n,i},v_{n,i})=-x_{k(n),i}.
\eneq
Let $v_{n,i}=e^{(n,i)},$ if ${\rm sp}(z_{n,i})\not=S^1.$ By
\ref{TRbot} and applying \ref{group}
\beq\label{asum-5++}
&&\tau({1\over{2\pi
i}}\log((z_{n,i}v_{n,i}^*w_{k(n)}^*\af(z_{n,i})w_{n+k}v_{n,i}))\\
&=&\tau({1\over{2\pi
i}}\log(v_{n,i}(z_{n,i}v_{n,i}^*w_{k(n)}^*\af(z_{n,i})w_{k(n)}v_{n,i}v_{n,i}^*))\\
&=&\tau({1\over{2\pi
i}}\log(v_{n,i}z_{n,i}v_{n,i}^*z_{n,i}^*z_{n,i}w_{k(n)}^*\af(z_{n,i})w_{k(n)}))\\
&=&\tau(\log({1\over{2\pi}}v_{n,i}z_{n,i}v_{n,i}^*z_{n,i}^*))+
\tau(\log({1\over{2\pi}}z_{n,i}w_{k(n)}^*\af(z_{n,i})w_{k(n)}))\\
&=& \text{bott}_1(z_{n,i},\, v_{n,i})+
\tau(\log({1\over{2\pi}}z_{n,i}w_{k(n)}^*\af(z_{n,i})w_{k(n)}))=0.
\eneq

We then obtain a unitary $u_n$ (with the form of
$\oplus_{i=1}^{R(n)}w_{k(n)}{\rm diag}(v_{n,i}, \cdots ,
v_{n,i})$) such that
\beq\label{asum-5+1}
{\rm ad}\, u_n|_{D_n}={\rm id}_{D_n}, \,\,\,\|{\rm ad}\, u_n\circ
\af(z_{n,i})-z_{n,i}\|<\ep_n/2\\\label{asum-5+2}
 \tau({1\over{2\pi
i}}\log(z_{n,i}u_n^*\af(z_{n,i}^*)u_n))=0.
\eneq
$n=1,2,....$ By (\ref{asum-4+2})and choice of $\dt_n'$, we may
also assume  that
\beq\label{asum-6+1}
\|{\rm ad}\, u_n(\af(a))-a\|<\ep_n/2\rforal a\in {\cal G}_n.
\eneq
By applying \ref{PAT}, we also have that
\beq\label{asum-6}
\tau({1\over{2\pi
i}}\log(z_{n,j}u_{n+1}^*\af(z_{n,j}^*)u_{n+1})=0.
\eneq



 Put
\beq\label{asum-10}
b_1&=&{1\over{2\pi i}}\log({\rm ad}\,
u_{n+1}^*(z_{n,i})\af(z_{n,i}^*))\\
b_2 &=&{1\over{2\pi i}}\log({\rm ad}\,
u_n^*(z_{n,i})\af(z_{n,i}^*)))\\
 b_3&=&{1\over{2\pi
i}}\log(z_{n,i}{\rm ad}\, u_{n}^*u_{n+1}(z_{n,i}^*)).
\eneq

Then, by (\ref{asum-5+2}) and (\ref{asum-6}),
\beq\label{asum-11}
\tau(b_1)
&=& \tau({1\over{2\pi i}}\log(z_{n,i}{\rm
ad}(u_{n+1})(\af(z_{n,i})))=0\andeqn\\
\tau(b_2)
&=& \tau(\log({1\over{2\pi i}}z_{n,i}{\rm ad}\,
u_n(\af(z_{n,i}^*))))=0.
\eneq
But
\beq\label{asum-11+1}
{\rm ad}\, u_{n+1}^*(e^{2\pi i h_3})=e^{2\pi i b_1}e^{-2\pi ib_2}.
\eneq
Therefore, by applying (\ref{group}),
\beq\label{asum-12}
\tau(h_3)&=& \tau({1\over{2\pi i}}\log( {\rm ad}\,
u_{n+1}^*(e^{2\pi i h_3}))=\tau({1\over{2\pi i}}\log(e^{2\pi i
b_1}e^{-2\pi ib_2}))\\
&=& {1\over{2\pi i}}\tau(2\pi i b_1)+{1\over{2\pi i}}\tau(2\pi ib_2)=0 .
\eneq
Note also, by (\ref{asum-6+1}) ,
\beq\label{asum-12-1}
&&\|[j(a),u_n^*u_{n+1}]\|<\ep_n/2+\ep_{n+1}/2\tforal a\in {\cal
G}_n
\eneq
 In particular,
 \beq
&&\|[z_{n,i},\,
u_n^*u_{n+1}]\|<\ep_n/2+\ep_{n+1}/2.
\eneq

However, by \ref{TRbot}
\beq\label{asum-13}
\tau(\text{bott}_1(u_{n+1}^*u_n,\,z_{n,i}))=\tau(h_3)=0.
\eneq
Since ${\rm ker}\rho_C=\{0\},$
$$
\text{bott}_1(u_{n+1}^*u_n,\,z_{n,i}))=0.
$$
 Therefore
\beq\label{asum-14}
\text{bott}_1(z_{n,i}, u_n^*u_{n+1})=0.
\eneq
Hence we may assume that
\beq\label{asum-15}
\text{bott}_1(j, u_n^*u_{n+1})|_{{\cal P}_n}=0.
\eneq
Consequently, (also using (\ref{asum-12-1})),  for
sufficiently large $n$ (so $\ep_n$ is sufficiently small), one may
apply \ref{BHPT} ( $K_i(B)$ is torsion free and $K_1(A)=\{0\}$ and $K_0(A)$ is divisible)
to obtain a continuous path of unitaries
$\{w_n(t): t\in [0,1]\}$ such that
\beq\label{asum-16}
&&w_n(0)=1,\,\,\, w_n(1)=u_n^*u_{n+1} \andeqn\\
&& \|[j(a),
w_n(t)]\|<1/n
\eneq
for all  $a\in {\cal F}_n$ and $t\in [0,1].$

Define
\beq\label{asum-17}
U(t+n)=u_nw_n(t)\tforal t\in [0,1].
\eneq
Then $U(t)$ is continuous on $[0, \infty).$ One then easily checks
that
$$
\lim_{n\to\infty}{\rm ad}\, U(t)\circ\af(a)=a\tforal a\in B.
$$

\end{proof}

\begin{thm}\label{ATT}
Let $B$ be a unital A$\T$-algebra and let $\af\in Aut(B).$ Then
$A\rtimes_{\af}\Z$ can be embedded into a unital simple AF-algebra
if and only if $A$ admits a faithful $\af$-invariant tracial
state.
\end{thm}

\begin{proof}
It suffices to show the ``if" part of the statement.  By
\ref{CembAH}, there is a unital simple AF-algebra $B$ with a
unique tracial state $\tau$ and a unital monomorphism $h: A\to B$
such that
$$
\tau\circ h\circ \af=\tau\circ h.
$$
It follows from \ref{Asmu} that there exists a unital simple
AF-algebra $C,$ a unital monomorphism $\phi: A\to C$ and a
continuous path of unitaries $\{U(t): t\in [0, \infty)$ such that
\beq\label{ATT1}
\lim_{t\to\infty}{\rm ad}\, U(t)\circ \phi\circ
\af(a)=\phi(a)\rforal a\in A.
\eneq
It follows from \ref{MC} that $A$ can be embedded into a unital
simple AF-algebra.

\end{proof}

\begin{cor}\label{CATT}
Let $A$ be a unital simple A$\T$-algebra and let $\af\in Aut(A).$ Then
$A\rtimes_{\af}\Z$ can always be embedded into a unital simple AF-algebra.
\end{cor}

\begin{proof}
Since $A$ is simple, $A$ always admits a faithful $\af$-invariant tracial state.
\end{proof}

\begin{rem}
{\rm This section is based on \cite{M1}. It simplifies  and
clarifies some arguments used in \cite{M1} for the purpose of
proving \ref{ATT}. It should be noted that a unital simple
A$\T$-algebra may not have real rank zero. Theorem \ref{ATT} also
applies to non-simple A$\T$-algebras.

}

\end{rem}

\section{Asymptotic embeddings  }

\begin{lem}\label{1emb}
Let $A$ be a unital \CA\, in ${\cal N},$ let $\af\in Aut(A)$ and
let $C$ be a unital separable simple AF-algebra. Let $T:
A_+\setminus \{0\}\to \N\times \R_+.$ Suppose that there is a
$T$-full unital asymptotic  embedding $\{\phi_n\}$ from $A$ to $C$
such that, for any finite subset ${\cal P}\subset K_0(A),$
\beq\label{1L-0}
[\phi_n]|_{\cal P}=[\phi_n\circ \af]|_{\cal P}
\eneq
for all sufficiently large $n.$ Then there is a unital  asymptotic
embedding $\{\psi_n\}$
 from $A$ to the unital simple  AF-algebra $B=C\otimes {\cal U}$
and there is a sequence of unitaries $u_n\in B$ such that
\beq\label{1L1}
\lim_{n\to\infty}\|\psi_n\circ \af(a)-{\rm ad}\, u_n\circ
\psi_n(a)\|=0
\eneq
for all $a\in A.$ If $C$ has a unique tracial state, one can
require that $B$ also has a unique tracial state. Moreover,
suppose that $A_0\subset A$ is a unital \SCA\, and for each finite
subset ${\cal Q}\in \underline{K}(A_0),$
\beq\label{1L1-1-1}
[\phi_n]|_{[j_0](\cal Q)}=[\phi_m]|_{[j_0](\cal Q)}
\eneq
for all $n, m\ge n_0$ for some $n_0>0,$ where $j_0: A_0\to A$ is
the embedding, then we may also require that
\beq\label{1L1-1}
[\psi_n]|_{[j_0](\cal Q)}=[\psi_m]|_{[j_0](\cal Q)}
\eneq
for all sufficiently large $n$ and $m.$

\end{lem}

\begin{proof}
We may assume that
$\phi_n(1_A)=1_B.$

Since $\{\phi_n\}$ is $T$-full, by \ref{fullP1} and \ref{inj1},
for any finite subset ${\cal G}_0\subset A,$ there exists $n_0>0$
such that
\beq\label{1emb1}
\|\phi_n(a)\|\ge (1/2)\|a\|\tforal a\in {\cal G}_0
\eneq
for all $n\ge n_0.$
 Fix a finite subset ${\cal F}\subset A$ and
$\ep>0.$ Let ${\cal G}\subset A,$ ${\cal P}\subset
\underline{K}(A)$ be finite subsets, $\dt>0$ and $N>0$ be integer
meet the requirements of \ref{Lunq}. Fix ${\cal Q}\subset
\underline{K}(A_0).$ Without loss of generality, we may assume
that $[j_0]({\cal Q})\subset {\cal P}.$

Let $K=N(N+1).$ Choose ${\cal G}'=\cup_{j=0}^N\af^j({\cal G}).$
Note that $\pi\circ \{\phi_n\circ \af^j\}: A\to q_{\infty}(\{C\})$
is full for each $j.$ Therefore there are $T_j:
A_+\setminus\{0\}\to (\N, \R_+),$ $j=0,1,2,...,$ such that
$\{\phi_n\circ \af^j\}$ are $T_j$-full. Let ${\cal S}_j={\cal
S}_j(\ep, {\cal F}, A, T_j)$ be a finite subset required by
\ref{Lunq}, $j=0, 1,2,...,N.$

Suppose that $\phi\subset \{\phi_n\}$  such that it is a $({\cal
G}',{\cal P},\dt/2)$-map. By (\ref{1emb1}), we may also assume
that
\beq\label{1emb2}
\|\phi(a)\|\ge (1/2)\|a\|\tforal a\in {\cal G}'.
\eneq

We may also assume that
\beq\label{Nt1}
[\phi\circ \af^j]|_{\cal P}=[\phi]|_{\cal P},\,\,\,j=1,2,...,N.
\eneq
Define
$$
N_a=\max_{0\le j\le K} T(\af^j(a))_\N\andeqn R_a=\max_{0\le j\le
K}T(\af^j(a))_\R.
$$
Define $T_1: A_+\setminus \{0\}\to \N\times \R_+$ by $T_1(a)=(N_a,
R_a)$ for $a\in A_+\setminus \{0\}.$ Thus we can
 choose $\phi\in \{\phi_n\}$ such that $\phi\circ
\af^j$ is $T_1$-${\cal S}_j$ full for $j=1,2,...,N.$

Write $C=C\otimes {\cal U}.$ Define $\psi: A\to C\otimes M_{K}$ by
\beq\label{Nt2}
\psi(a)={\rm diag}(\Phi(a),\Phi(\af(a)),...,\Phi(\af^N(a))),
\eneq
for $a\in A,$ where
\beq\label{Nt3}
\Phi(b)=\phi(b)\otimes 1_{M_N}={\rm
diag}(\overbrace{\phi(a),\phi(a),...,\phi(a)}^N)
\eneq
for $b\in A.$ Thus
\beq\label{Nt4}
\psi\circ \af(a)={\rm
diag}(\Phi(\af(a)),\Phi(\af^2(a)),...,\Phi(\af^{N+1}(a)))
\eneq
for $a\in A.$ Note that there is a permutation matrix $W\in
M_{N+1}$ such that
\beq\label{Nt5}
{\rm ad}\, W\circ \psi(a)={\rm diag}(\Phi(\af(a)),
\Phi(\af^2(a)),...,\Phi(\af^N(a)), \Phi(a))
\eneq
for all $a\in A.$

Since each $\phi\circ \af^j$ is $T_1$-${\cal S}_j$ full, by
counting multiplicity of each $\phi(\af^j(a))$ for $j=1,2,...,N,$
using (\ref{Nt1}) and by applying \ref{Lunq}, we obtain a unitary
$U\in C\otimes M_{K}$ such that
\beq\label{Nt6}
{\rm ad}\, U\circ \psi\approx_{\ep} \psi\circ \af \,\,\,{\rm
on}\,\,\,{\cal F}
\eneq
By embedding  $M_K$ unitally into ${\cal U},$ we obtain a $({\cal G}, {\cal
P}, \dt)$-map from $A$ into $C\otimes {\cal U}.$

By (\ref{Nt1}) and (\ref{Nt2}), we see that, for any $x\in {\cal
P},$
\beq\label{Nt7}
[\psi](x)=[\phi](x)\otimes 1_{\cal U}.
\eneq
This proves (\ref{1L1-1}). Moreover,
$$
\|\psi(a)\|\ge (1/2)\|a\|\tforal a\in {\cal F}.
$$


\end{proof}

\begin{cor}\label{C1emb}
If in addition to the hypothesis in \ref{1emb}, there are
$\af_1, \af_2, ..., \af_k\in Aut(A)$ such that
$\af_j\circ \af_i=\af_i\circ \af_j$ ($i,j=1,2,...,k$) and
\beq\label{C1emb-1}
\tau\circ \phi_1(a)=\tau\circ \phi_n(a)=\tau\circ \phi_n\circ \af_j(a)
\eneq
for all $a\in A$ and for all $\tau\in T(C),$ $n=1,2,...$ and $j=1,2,...,k,$
where $t$ is the unique tracial state of ${\cal U},$
then  we may also require (in addition to the conclusion of \ref{1emb}) that
\beq\label{C1emb-2}
(\tau\otimes t)\circ \psi_1(a)=(\tau\otimes t)\circ \psi_n(a)=(\tau\otimes t)\circ \psi_n\circ \af_j(a)
\eneq
for all $a\in A$ and for all $\tau\in T(C),$ $n=1,2,...$ and $j=1,2,...,k.$
\end{cor}

\begin{proof}
Note that, by (\ref{Nt2}),
\beq\label{C1emb-3}
\psi\circ \af_j(a)={\rm diag}(\Phi(\af_j), \Phi(\af\circ \af_j(a)),\cdots, \Phi(\af^N\circ \af_j(a))).
\eneq
By (\ref{C1emb-1}) (and (\ref{Nt3})), we have
\beq\label{C1emb-4}
(\tau\otimes t')\circ \Phi(\af^i\circ \af_j(a)))=(\tau\otimes t')\circ \Phi(\af_j(a))=(\tau\otimes t')\circ \Phi(a)
\eneq
for all $\tau\in  T(C),$ where $t'$ is the tracial state on $M_K.$
It is then easy  to  conclude that (\ref{C1emb-2}) holds.

\end{proof}

\begin{lem}\label{2L1}
Let $A$ be a unital separable amenable \CA, let $\af\in Aut(A)$
and let $B$ be a unital \CA. Suppose that there exists a sequence
of \morp s $\{\phi_n\}$ from $A$ to $B$ which satisfies the
following:

{\rm (1)} $\lim_{n\to\infty}\|\phi_n(ab)-\phi_n(a)\phi_n(b)\|=0$
for all $a, b\in A,$

{\rm (2)} there exists a sequence of unitaries $\{v_n\}$ in $B$
such that
\beq\label{2L-1}
\lim_{n\to\infty}\|{\rm ad}\, v_n\circ
\phi_n(a)-\phi_n(\af(a))\|=0
\eneq
for all $a\in A.$

Then there exists a sequence of \morp s $\{\psi_n\}$ from
$A\rtimes_\af\Z$ to $B$ such that
\beq\label{2L-2}
\lim_{n\to\infty}\|\psi_n(xy)-\psi_n(x)\psi_n(y)\|=0
\eneq
for all $x,y\in A\rtimes_\af\Z,$
\beq\label{2L-3}
\lim_{n\to\infty}\|\psi_n(a)-\phi_n(a)\|=0\andeqn
\lim_{n\to\infty}\|\psi_n(u_\af)-v_n\|=0
\eneq
for all $a\in A.$
\end{lem}

\begin{proof}
Define
\beq\label{2L-4}
L_n(\sum_{-N\le k\le N} a_ku_\af^k)=\sum_{-N\le k\le
N}\phi_n(a_k)v_n^k
\eneq
for $a_k\in A.$ For any $a, b\in A$ and any integers $m$ and $k,$
\beq\label{2L-5}
L_n(au_\af^mbu_\af^k)=L_n(a\af^{-m}(b)u_\af^{k+m})=\phi_n(a\af^{-m}(b))v^{k+m}
\andeqn
\eneq
\beq\label{2L-5+}
L_n(au_\af^m)L_n(bu_\af^k)&=&\phi_n(a)v_n^m\phi_n(b)v_n^k\\
&=&[\phi_n(a)v_n^m\phi_n(b)v^{-m}]v_n^{k+m}\
\eneq
However, by (\ref{2L-1}) and (\ref{2L-2}),
\beq\label{2L-6}
\lim_{n\to\infty}\|\phi_n(a\af^{-m}(b))v_n^{k+m}-[\phi_n(a)v_n^m\phi_n(b)v^{-m}]v_n^{k+m}\|=0
\eneq
It follows that
\beq\label{2L-7}
\lim_{n\to\infty}\|L_n(au_\af^mbu_\af^k)-L_n(au_\af)L_n(bu_\af^k)\|=0
\eneq
Therefore, for any $x, y\in A\rtimes_\af\Z,$
\beq\label{2L-8}
\lim_{n\to\infty}\|L_n(xy)-L_n(x)L_n(y)\|=0.
\eneq
Define $L: A\rtimes_{\af}\Z\to l^{\infty}(\{B\})$ by
$L(c)=\{L_n(c)\}$ for $c\in A\rtimes_{\af}\Z.$ Then $\pi\circ L:
A\rtimes_{\af}\Z\to q_{\infty}(\{B\})$ is a \hm. Since $A$ is
assumed to be amenable, so is $A\rtimes_\af\Z.$ It follows from
\cite{CE} that  there is a sequence of \morp s $\psi_n:
A\rtimes_{\af}\Z\to B$ such that
\beq\label{2L-9}
\lim_{n\to\infty}\|L_n(x)-\psi_n(x)\|=0
\eneq
for all $x\in A\rtimes_{\af}\Z\to B.$  This proves the lemma.
\end{proof}

\begin{cor}\label{1C}
Let $A\in {\cal N}$ be a unital \CA, let $\af\in Aut(A)$ and let
$C$ be a unital simple AF-algebra. Let $T: A_+\setminus\{0\}\to
(\N, \R_+)$ be a map. Suppose that there is a $T$-full unital
asymptotic  embedding $\{\phi_n\}$ from $A$ to $C$ such that, for
any finite subset ${\cal P}\subset K_0(A),$
$$
[\phi_n]|_{\cal P}=[\phi_n\circ \af]|_{\cal P}
$$
for all sufficiently large $n.$ Then  $A\rtimes_{\af}\Z$ can be
unitally asymptotically embedded into a unital simple AF-algebra
\end{cor}

\begin{proof}
By embedding $C$ into $C\otimes {\cal U},$ we may assume that $C\cong
C\otimes {\cal U}.$ In particular, $K_0(C)$ is divisible. The corollary
then follows from \ref{1emb} and \ref{2L1}.

\end{proof}

\begin{cor}\label{1C2}
Let $A\in {\cal N}$ be a unital \CA, let $\af\in Aut(A)$ and let
$C$ be a unital simple AF-algebra. Suppose that there is a
monomorphism $\phi: A\to C$
$$
\phi_{*0}=(\phi\circ \af)_{*0}.
$$
Then there exists a unital simple AF-algebra $B,$ a sequence of
unital monomorphisms $\{\phi_n\}$ and a sequence of unitaries
$\{u_n\}$ in $B$ such that
$$
\lim_{n\to\infty}\|{\rm ad}\, u_n \circ \phi_n(a)-\phi_n\circ
\af(a)\|=0
$$
for all $a\in A.$

Moreover, if  $C$ has a unique tracial state $\tau$ such that
$$
\tau\circ \phi(a)=\tau\circ \phi(\af(a))
$$
for all $a\in A,$  then
there is an asymptotic embedding $\{\psi_n\}: A\rtimes_{\af}\Z\to B$ such that
$$
\lim_{n\to\infty}\|\psi_n\circ j(a)-\phi_n(a)\|=0\tforal a\in A.
$$
We may assume that $B$ has a unique tracial
state $t$ and
$$
\lim_{n\to\infty}t(\phi_n(a))=\lim_{n\to\infty}t(\phi_n(\af(a)))=\tau(\phi(a))
$$
for all $a\in A$ as well as
\beq\label{1C2+1}
[\phi_n]=[\phi_n\circ \af]\,\,\, in \,\,\, KL(A, B).
\eneq
\end{cor}

\begin{proof}
The  proof is evident by replacing $\phi_n$ in the proof of
\ref{1emb} by the monomorphism $\phi.$ We also apply \ref{C1emb}.

\end{proof}

\begin{prop}\label{T-full}
Let $A$ be a unital separable amenable \CA\, and let $B$ be a
unital \CA. Suppose that $\af\in Aut(A)$ satisfies the cyclic
Rokhlin property (\ref{DU}). Suppose that there is a sequence
$\{\phi_n\}$ of \morp s from $A\rtimes_{\af} \Z $ to $B$ for which
$$
\lim_{n\to\infty}\|\phi_n(ab)-\phi_n(a)\phi_n(b)\|=0\tforal a,b\in
A.
$$
 Then $\{\phi_n\}$ is
$T$-full for some $T: (A\rtimes_{\af}\Z)_+\setminus \{0\}\to (\N,
\R_+)$ if and only if $\{\phi_n\circ j\}$ is $T'$-full for some
$T': A_+\setminus \{0\}\to  (\N, \R_+),$ where $j: A\to
A\rtimes_\af\Z$ is the natural embedding.
\end{prop}

\begin{proof}
It suffices to show the ``if" part. Suppose that, for some $T':
A_+\setminus \{0\}\to  (\N, \R_+),$ $\{\phi_n\circ j\}$ is $T'$-
full. Let $\pi: l^{\infty}(\{B\})\to q_{\infty}(\{B\})$ be the
quotient map. Put ${\bar \Phi}=\pi\circ \{\phi_n\}.$ Then ${\bar
\Phi}: A\rtimes_{\af}\Z\to q_{\infty}(\{B\})$ is a \hm. It follows
from (\ref{fullP1}) that ${\bar \Phi}\circ j$ is full. Let $J_1$
be a closed two-sided ideal of $q_{\infty}(\{B\}).$ Suppose that
$$
J=\{b\in A\rtimes_\af\Z: {\bar \Phi}(b)\in J_1\}.
$$
Let $I=J\cap j(A).$ Since $\af$ is cyclic Rokhlin, by a result of
Pasnicu and Phillips (see 4.1 of \cite{Lnemb2}), $J$ is proper  if
and only if $I$ is proper. Since ${\bar \Phi}|_A$ is a full,
$I=\{0\}.$ So is $J.$ This implies that, for any $b\not=0,$ ${\bar
\Phi}(b)$ does not belong to any closed two-sided ideal of
$q_{\infty}(\{B\}).$ Thus ${\bar \Phi}$ is full. Therefore, by
\ref{fullP1}, there is $T: (A\rtimes_{\af}\Z)_+\setminus \{0\}\to
(\N, \R_+)$ such that $\{\phi_n\}$ is $T$-full.

\end{proof}

\begin{cor}\label{CT-full}
Let $A$ be a unital separable amenable \CA\, and let $B$ be a
unital \CA. Suppose that $\af\in Aut(A)$ satisfies the cyclic
Rokhlin property (\ref{DU}). Suppose that there is an asymptotic
embedding $\{\phi_n\}: A\rtimes_{\af}\Z\to B$  and there is a full
monomorphism $h: A\to B$ such that
\beq\label{Tf1}
\lim_{n\to\infty}\|h(a)-\phi_n\circ j(a)\|=0
\eneq
for all $a\in A,$ where $j: A\to A\rtimes_\af\Z$ is the natural
embedding. Then $\{\phi_n\}$ is $T$-full for some  $T:
(A\rtimes_{\af}\Z)_+\setminus \{0\}\to (\N, \R_+).$
\end{cor}

\begin{proof}
It follows from (\ref{Tf1}) that ${\bar \Phi}\circ j=h.$ Thus
${\bar \Phi}\circ j$ is full.
\end{proof}

\section{Asymptotic Embeddings for $C(X)\rtimes_{\af}\Z^2$ and quasidiagonality}

\

\begin{thm}\label{TembX}
Let $X$ be a compact metric space and let $\Lambda: \Z^2\to Aut(C(X))$ be  a $\Z^2$
action on $X.$

Then the following are equivalent:

{\rm (1)} there exists a strictly positive $\Lambda$-invariant
Borel probability measure $\mu$ on $X,$

{\rm (2)} there exists a unital simple AF-algebra $B,$ there
exists a monomorphism $h: C(X)\to B$ and there exists a $T$-full
asymptotic embedding
 $\{\phi_n\}: C(X)\rtimes_{\Lambda} \Z^2\to B$ for some
 $T:A_+\setminus\{0\}\to (\N, \R_+)$
such that
\beq\label{aembX1}
\lim_{n\to\infty}\|h(f)-\phi_n\circ j(f)\|=0
\eneq
for all $f\in C(X),$ where $j: C(X)\to C(X)\rtimes_{\Lambda} \Z^2$
is the natural embedding.

\end{thm}

\begin{proof}

Suppose that (2) holds. Without loss of generality, we may assume
that $\phi_n(1_{C(X)\rtimes_{\Lambda}\Z^2})= 1_B.$ Let $\tau\in
T(C)$ be a tracial state. Then $\tau\circ h$ is a strictly
positive Borel probability measure. Since
$C(X)\rtimes_{\Lambda}\Z^2$ is separable, there exists a
subsequence $\{n_k\}$ and a tracial state $t$ of
$C(X)\rtimes_{\Lambda}\Z^2$ such that
$$
t(a)=\lim_{k\to\infty} \tau\circ \phi_{n_k}(a)
$$
for all $a\in C(X)\rtimes_{\Lambda} \Z^2.$ Therefore $t\circ j$ is
$\Lambda$-invariant. However, by (\ref{aembX1})
$$
t\circ j=\tau\circ j.
$$
Thus $\tau\circ j$ is a strictly positive $\Lambda$-invariant Borel
probability measure on $X.$

Now suppose that (1) holds.

It follows from Cor 3.3 and  Theorem 4.3 of \cite{Lnemb2} (and its proof) that there is a unital
simple AF-algebra $B$ with a unique tracial state $\tau$ and a
unital embedding $\phi: C(X)\otimes {\cal
U}\rtimes_{\gamma_1}\Z\to B,$
such that
$$
\tau\circ \phi\circ \imath (a)=\tau\circ \phi\circ \imath\circ
\gamma_j(a), \,j=1,2
$$
for all $a\in C(X)\otimes {\cal U},$ where $\gamma_j=\af_j\otimes \sigma$ and $\sigma$ is
as in \ref{DU}, $j=1,2$  and
$\imath:
C(X)\otimes{\cal U}\to C(X)\otimes {\cal U}\rtimes_{\gamma_1}\Z$
is the natural embedding. Let $\D=\rho_B(K_0(B)).$ There is a unital simple
AF-algebra $B_1$ with $(K_0(B_1), K_0(B_1)_+, [1_{B_1}])=(\D, \D_+, 1)$ and there is a unital monomorphism
$h_{00}: B\to B_1$ which maps $K_0(B)$ to $\D$ by the map $\rho_B.$ Thus we may assume
that $B=B_1,$ without loss of generality. In particular, now we assume that
${\rm ker}\rho_B=\{0\}.$
By replacing $B$ by $B\otimes {\cal U},$ we may also assume
that $K_0(B)$ is divisible.

Put $A_0=C(X)\otimes {\cal U}.$ Suppose that $e\in
M_k(A_0\rtimes_\bt Z)$ is a projection. Let $E:
A_0\rtimes_{\gamma_1} \Z \to A_0$ be the canonical conditional
expectation. Since $\gamma_1$ has cyclic Rokhlin property, by Cor 4.2  of \cite{Lnemb2}, there
is a tracial state $t$ of $A_0$ such that
\beq\label{embX-1}
t(E(a))=\tau\circ \phi(a)
\eneq
for all $a\in A_0\rtimes_{\gamma_1} \Z.$ Define ${\bar \gamma}_2:
A_0\rtimes_{\gamma_1}\Z\to A_0\rtimes_{\gamma_1}\Z$ by ${\bar
\gamma}_2(a)=\gamma_2(a)$ for $a\in A_0$ and ${\bar
\gamma_2}(u_{\gamma_1})=u_{\gamma_1},$ where $u_{\gamma_1}$ is the
canonical unitary in $A_0\rtimes_{\gamma_1}\Z$ implementing
$\gamma_1.$ Since $\gamma_1$ commutes with $\gamma_2,$ ${\bar
\gamma}_2$ defines an automorphism on $A_0\rtimes_{\gamma_1}\Z.$
For the canonical expectation $E: A_0\rtimes_{\gamma_1} \Z \to
A_0,$ one computes that
\beq\label{embX+1}
E({\bar\gamma_2}(a))=\gamma_2(E(a))
\eneq
for all $a\in A_0\rtimes_{\gamma_1} \Z .$

Since $\tau\circ \phi\circ \imath$ is $\gamma_1$ and $\gamma_2$
invariant, by applying (\ref{embX-1}) and (\ref{embX+1})), we
obtain that
\beq\label{embX2}
t(E({\bar
\gamma}_2(a)))=t(\gamma_2(E(a)))=t(E(a))\,\,\rforal\,\,\, a\in
A_0\rtimes_{\gamma_1}\Z.
\eneq

Now let $e\in M_k(A_0\rtimes_{\bt_1}\Z)$ be a projection. It
follows from (\ref{embX2}) that
\beq\label{embX3}
t(E({\bar \gamma}_2(e)))=t(E(e)).
\eneq
It follows that
\beq\label{embX4}
\tau\circ \phi({\bar \gamma_2}(e))=\tau\circ \phi(e).
\eneq
Since ${\rm ker}\rho_B=\{0\},$ we conclude that
\beq\label{embX5}
\phi_{*0}=(\phi\circ {\bar \gamma_2})_{*0}
\eneq
Since $K_1(B)=0$ and $K_0(B)$ is divisible, we have
\beq\label{ambX5-1}
[\phi]=[\phi\circ {\bar \gamma_2}]\,\,\,\text{in}\,\,\,KL(A_0\rtimes_{{\bar \gamma_1}}\Z, B).
\eneq
It follows from \ref{1emb},  \ref{2L1} and \ref{T-full}  that there exists  a $T$-full  (for some
$T: A_+\setminus\{0\}\to \N\times \R_+$)
asymptotic embedding $\{\psi_n\}$ from
$(A_0\rtimes_{\gamma_1}\Z)\rtimes_{{\bar \gamma}_2}\Z$  into $B\otimes {\cal U} (\cong B)$
 satisfying the
following.
\beq\label{embX6}
\lim_{n\to\infty}\|\psi_n(a)-h_n(a)\|=0\rforal a\in A_0\rtimes_{{\gamma}_1}\Z,
\eneq
where $h_n: A_0\rtimes_{\gamma_1}\Z\to B\otimes {\cal U}$ is a unital monomorphism.
Moreover,  by \ref{C1emb}, for each $n,$
\beq\label{embX6+}
(\tau\otimes t)\circ h_n(a)=(\tau\otimes t)\circ h_1(a)\rforal a\in A_0\rtimes_{{\gamma}_1}\Z.
\eneq
Again,
\beq\label{embX7}
[h_n]=[h_1]\,\,\,\text{in}\,\,\, KL(A_0\rtimes_{\gamma_1}\Z, B\otimes {\cal U}),
\eneq
$n=1,2,....$
It follows from (\ref{embX7}), (\ref{embX6+}) and \ref{cdT1} that there exists a sequence of unitaries
$\{w_n\}\subset B\otimes {\cal U}$ such that
\beq\label{embX8}
\lim_{n\to\infty}\|{\rm ad}\, w_n\circ h_n\circ j(f)-h_1\circ j(f)\|=0\rforal f\in C(X).
\eneq


Thus, (by also (\ref{embX6}),
\beq\label{embX9}
\lim_{n\to\infty}\|{\rm ad}\, w_n\circ \psi_n\circ
j_1(f)-\psi_1\circ j_1(f)\|=0
\eneq
for all $f\in C(X).$ Now by replacing $\psi_n$ by ${\rm ad}\,
w_n\circ \psi_n,$ we may assume that
\beq\label{embX10}
\lim_{n\to\infty}\| \psi_n\circ j(f)-h_1\circ j(f)\|=0\rforal f\in C(X).
\eneq

\end{proof}

\begin{cor}\label{Cqd}
Let $X$ be a compact metric space and let $\Lambda: \Z^2\to Aut(C(X))$ be a $\Z^2$-action
on $X.$  Suppose that there is a strictly positive $\Lambda$-invariant Borel probability measure
$\mu.$
Then $C(X)\rtimes_{\Lambda}\Z^2$ is quasidiagonal.
\end{cor}

We would also mention the following two asymptotic embeddings:

\begin{prop}\label{AHasym}
Let $A$ be a unital AH-algebra and let $\af\in Aut(A)$ be an
automorphism. Then there is a $T$-full (for some $T: A_+\setminus
\{0\}\to \N\times \R$) asymptotic embedding
$\{\phi_n\}:A\rtimes_{\af}\Z\to C$ for some unital simple
AF-algebra $C$ if and only if $A$ admits a faithful
$\af$-invariant tracial state.
\end{prop}

\begin{proof}
By replacing $A$ by $A\otimes {\cal U}$ and replacing $\af$ by $\af\otimes \sigma,$
we may assume that $\af$ has the cyclic Rokhlin property of \ref{DU}.
By \ref{CembAH}, there is a unital monomorphism $h: A\to B$ for some unital simple
AF-algebra $B$ with unique tracial state $\tau$ such that ${\rm ker}\rho_{B}=\{0\}$ such that
\beq\label{AHasym1}
\tau\circ h=\tau\circ h\circ \af
\eneq
for all $\tau\in T(B).$ Since ${\rm ker}\rho_B=\{0\},$ it follows that
\beq\label{AUasym2}
h_{*0}=(h\circ\af)_{*0}
\eneq
It follows from \ref{1C2} that there is an asymptotic embedding
$\{\phi_n\}: A\rtimes_{\af}\Z\to B\otimes {\cal U}$ and a sequence of monomorphisms $h_n: A\to B\otimes {\cal U}$
such that
$$
\lim_{n\to\infty}\|\phi_n\circ j(a)-h_n(a)\|=0\andeqn
$$
$$
\tau\circ h_n(a)=\tau\circ h_1(a)\tforal a\in A.
$$
Note since $K_1(B)=\{0\}$ and $K_0(B\otimes {\cal U})$ is divisible and ${\rm ker}\rho_B=\{0\},$
$$
[h_n]=[h_1]\,\,\,\text{in}\,\,\,KL(A, B\otimes {\cal U}).
$$
It follows from \ref{cdT1} that there is a sequence of unitaries $\{w_n\}\subset B$ such that
$$
\lim_{n\to\infty}\|{\rm ad}\, w_n\circ h_n(a)-h_1(a)\|=0\tforal a\in A.
$$
We now replace $\psi_n$ by ${\rm ad}\, w_n\circ \psi_n.$ Thus, we may assume that
$$
\lim_{n\to\infty}\|\psi_n\circ j(a)-h_1(a)\|=0\tforal a\in A.
$$
Since $B$ is simple, it is easy to see that there is $T': A_+\setminus\{0\}\to \N\times \R_+$ such that
$\{\psi_n\circ j\}$ is $T'$-full. Since we have assume that $\af$ has the cyclic Rokhlin property,
by \ref{T-full}, $\{\psi_n\}$ is $T$-full for some $T: A_+\setminus\{0\}\to \N\times \R_+.$

\end{proof}

\begin{prop}\label{Res}
Let $A$ be a separable amenable residually finite dimensional \CA\, in ${\cal N}$ and $\af\in Aut(A)$ such that
$(\af)_{*0}={\rm id}|_{K_0(A)}.$ Then there is an asymptotic embedding $\{\phi_n\}: A\rtimes_{\af}\Z
\to B$ for some unital simple AF-algebra.
\end{prop}

\begin{proof}
It follows from \cite{Lnemb1} that there is a unital monomorphism $h: A\to B$ for some unital simple AF-algebra
$B$ with a unique tracial state $\tau.$
By mapping $B$ to another simple AF-algebra, if necessary, we may assume that $K_0(B)$ is divisible.
Since $(\af)_{*0}={\rm id}|_{K_0(A)},$ it follows that
\beq\label{Res1}
[h\circ \af]=[h]\,\,\,\text{in}\,\,\,KL(A,B).
\eneq
It follows from \ref{1C} that $A\rtimes_{\af}\Z$ can be asymptotically embedded into a unital simple AF-algebra $B.$

\end{proof}

\vspace{0.4in}

\noindent
Department of Mathematics\\
East China Normal University\\
Shanghai, China\\
and (current)\\
Department of Mathematics\\
University of Oregon\\
Eugene, Oregon 97403\\
U.S.A.


\begin{thebibliography}{BH}

\bibitem{BK1} B.  Blackadar and E.  Kirchberg {\em Generalized inductive limits of finite-dimensional $C\sp *$-algebras},  Math. Ann. {\bf 307} (1997), 343--380.


\bibitem{BK2} B.  Blackadar and E.  Kirchberg {\em Inner quasidiagonality and strong NF algebras}, Pacific J. Math. {\bf 198} (2001), 307--329.

\bibitem{BlD} L. G.\,ÊBrown and M. D\u{a}d\u{a}rlat, {\em Extensions of $C\sp *$-algebras and quasidiagonality},  J. London Math. Soc. {\bf  53} (1996), 582--600.

\bibitem{Bn1} N. Brown, {\em AF embeddability of crossed products of AF algebras by the integers},  J. Funct. Anal.  {\bf 160}  (1998),
  no. 1, 150--175.

\bibitem{Bn3} N. Brown, {\em Crossed products of UHF algebras by some amenable groups}, Hokkaido Math. J.
 {\bf 29} (2000), 201--211


\bibitem{Bn2} N. Brown, {\em On quasidiagonal $C\sp *$-algebras.  Operator algebras and applications}
  19--64, Adv. Stud. Pure
Math., 38, Math. Soc. Japan, Tokyo, 2004.


\bibitem{BnD} N. Brown and M. Dadarlat, {\em Extensions of quasidiagonal $C\sp *$-algebras and $K$-theory}, Operator algebras and applications, 65--84, Adv. Stud. Pure Math., 38, Math. Soc. Japan, Tokyo, 2004.

\bibitem{CE} M.D. Choi and E. Effros, {\em The completely positive lifting problem for \CA s}, Ann. of Math.
{\bf 104} (1976), 309-322.

\bibitem{D1} M. D\u{a}d\u{a}rlat, {\em Quasidiagonal morphisms and homotopy},  J. Funct. Anal. 151 (1997),  213--233.

\bibitem{D2} M. D\u{a}d\u{a}rlat, {\em On the approximation of quasidiagonal $C\sp *$-algebras},  J. Funct. Anal. {\bf 167 }(1999), 69--78.

\bibitem{D3} M. D\u{a}d\u{a}rlat, {\em Residually finite dimensional $C\sp *$-algebras and subquotients of the CAR algebra},  Math. Res. Lett. {\bf 8} (2001),  545--555.



\bibitem{DHS} K. Davidson, D.  Herrero and N.  Salinas, {\em Quasidiagonal operators, approximation, and $C\sp *$-algebras},   Indiana Univ. Math. J. {\bf 38} (1989), 973--998.

\bibitem{ELP} S.  Eilers, T.  Loring and G. K. Pedersen, {\em Quasidiagonal extensions and $AF$ algebras},  Math. Ann. {\bf 311} (1998),  233--249.

\bibitem{E1}G. A.  Elliott, {\em On the classification of $C\sp *$-algebras of real rank zero},
 J. Reine Angew. Math. {\bf 443} (1993), 179--219.

\bibitem{EG} G. A. Elliott and G.  Gong, {\em On the classification of $C\sp *$-algebras of real rank zero. II},
 Ann. of Math. {\bf 144} (1996), 497--610

\bibitem{Ex}R. Exel, {\em The soft torus and applications to
almots commuting matrics}, Pacific J. Math., {\bf 160} (1993),
207-217.


\bibitem{GL} G.  Gong and H. Lin, {\em Almost multiplicative morphisms and $K$-theory},
 Internat. J. Math. {\bf 11} (2000),  983--1000.


\bibitem{Hd} D.  Hadwin, {\em Strongly quasidiagonal $C\sp *$-algebras},   J. Operator Theory  {\bf 18}  (1987),
 3--18.

\bibitem{hlx} S. Hu, H.  Lin and Y.  Xue, {\em The tracial topological rank of $C\sp *$-algebras. II},
 Indiana Univ. Math. J. {\bf 53} (2004),  1578--1603.

\bibitem{K} A. Kishimoto, {\em Automorphisms of ${\rm A}\bold T$ algebras with the Rohlin property}, J. Operator Theory {\bf 40} (1998), 277--294.

\bibitem{KK} A.  Kishimoto and A.  Kumjian, {\em The Ext class of an approximately inner automorphism. II}, J. Operator Theory {\bf 46} (2001),  99--122.

\bibitem{Li} L. Li, {\em Classification of simple $C\sp *$-algebras: inductive limits of matrix algebras over one-dimensional spaces},  J. Funct. Anal. {\bf 192} (2002),  1--51.

\bibitem{Lnplms} H. Lin, {\em Tracial topological ranks of \CA s},
Proc. London Math. Soc., {\bf 83} (2001), 199-234.

\bibitem{Lnbk} H. Lin, {\em An introduction to the classification of amenable $C\sp *$-algebras}, World Scientific Publishing Co., Inc., River Edge, NJ, 2001. xii+320 pp. ISBN: 981-02-4680-3.

\bibitem{Lnemb1} H. Lin, {\em Residually finite dimensional and AF-embeddable $C\sp *$-algebras},  Proc. Amer. Math. Soc. {\bf 129} (2001), 1689--1696

\bibitem{Lnann} H. Lin, {\em Classification of simple $C\sp *$-algebras and higher dimensional
noncommutative tori},   Ann. of Math. (2) 157 (2003), no. 2, 521--544.



\bibitem{Lnduke} H. Lin,
          {\em Classification of simple $C^*$-algebras with
tracial topological rank zero},
  Duke Math. J.,{\bf 125} (2004), 91-119.

\bibitem{Lntr0} H. Lin, {\em Traces and simple $C\sp *$-algebras with tracial topological rank zero},
 J. Reine Angew. Math. 568 (2004), 99--137.


\bibitem{Lnuct} H.  Lin, {\em An approximate universal coefficient theorem},  Trans. Amer. Math. Soc. {\bf 357} (2005),  3375--3405.

\bibitem{Lngpots} H. Lin, {\em The Rokhlin property for automorphisms on a unital simple \CA }, Contemp.
Math., {\bf 414}, 2006, 189-215.

\bibitem{Lncd} H. Lin {\em Classification of \hm s and dynamical systems}, Trans. Amer. Math. Soc. {\bf 359} (2007), 859-895.

\bibitem{Lnemb2} H. Lin, {\em Embedding crossed products into a unital simple AF-algebra},
preprint, arxiv,org/ OA/0604047

\bibitem{Lnhomp} H. Lin, {\em  Approximate homotopy of \hm s from $C(X)$ into a simple \CA\,},
preprint, arxiv.org/ OA/0612125.

\bibitem{M1} H. Matui, {\em AF embeddability of crossed products of AT algebras by the integers and its application},
  J. Funct. Anal. {\bf 192} (2002),  562--580.

\bibitem{PcP} C. Pasnicu and N. C. Phillips {\em  Crossed products by Z with Rokhlin property }, in preparation.

\bibitem{P1} N. C. Phillips, {\em Crossed products of the Cantor set by free minimal actions of $\Bbb Z\sp d$},
  Comm. Math. Phys. {\bf 256} (2005), 1--42.

\bibitem{Pi} M. Pimsner, {\em  Embedding some transformation group $C\sp{*} $-algebras into AF-algebras},   Ergodic
Theory Dynam. Systems  {\bf 3}  (1983),  613--626.


 \bibitem{R1} M. R\o rdam, {\em On the structure of simple $C\sp *$-algebras tensored with a UHF-algebra. II}  J. Funct. Anal. {\bf 107} (1992),  255--269.

\bibitem{Ros} J. Rosenberg, {\em Appendix to Strongly quasidiagonal \CA s}, J. Operator Theory {\bf 18} (1987),
3-18.
Ê
\bibitem{Sa} N. Salinas, {\em Relative quasidiagonality and $KK$-theory}, Houston J. Math., {\bf 18} (1992), 97-116.

\bibitem{Th} F. J. Thayer, {\em Quasidiagonal \CA s}, J. Funct. Anal. {\bf 25} (1977), 50-57.

 \bibitem{Vj} J. Villadsen, {\em On the stable rank of simple $C\sp *$-algebras}, J. Amer. Math. Soc. {\bf 12} (1999) 1091--1102.

\bibitem{V1} D. Voiculescu, {\em A non-commutative Wely-con Nuemann theorem}, Rev. Roumanine Math.,
Pures Apply. {\bf 21} (1976), 97-113.

\bibitem{V2} D. Voiculescu, {\em Almost inductive limit automorphisms and embeddings into AF-algebras}, Ergodic
Theory Dynam. Systems, {\bf 6} (1986), 475-484.

\bibitem{V3} D. Voiculescu, {\em A note on quasi-diagonal $C\sp *$-algebras and homotopy},
  Duke Math. J.  {\bf 62}  (1991),   267--271.

\bibitem{V4} D. Voiculescu, {\em Around quasidiagonal operators}  Integral Equations Operator Theory  {\bf 17}
  (1993), 137--149.


 \bibitem{Z} S.  Zhang, {\em $K\sb 1$-groups, quasidiagonality, and interpolation by multiplier projections}, Trans. Amer. Math. Soc. {\bf 325} (1991),  793--818.


\end{thebibliography}
\end{document}